\documentclass[twoside]{siamltex}

\usepackage{amssymb}
\usepackage{amsmath}
\usepackage[final]{graphicx}
\usepackage{psfrag}
\usepackage{epsfig}
\usepackage{latexsym}
\usepackage{exscale}

\font\Bbb=msbm10 %at 12pt

\newcommand{\Rset}{\mbox{\Bbb R}}

\newcommand{\weg}[1]{ }

\def \dt    { h }
\def \dx    { {{\Delta} x} }

\def \1{1\!{\rm l}}
\def \A{{\cal O}\!\!\iota}
\def \half{{\textstyle\frac12}}
\def \z{\zeta}

\title{Fast and oblivious convolution quadrature}

\author{Achim Sch\"adle\footnotemark[1] \and
Mar\'ia L\'opez-Fern\'andez\footnotemark[5] \and
Christian Lubich\footnotemark[7]
}

\begin{document}

\maketitle

\renewcommand{\thefootnote}{\fnsymbol{footnote}}
\footnotetext[1]{ZIB Berlin, Takustr.~7, D-14195 Berlin, Germany.
~E-mail: {\tt schaedle@zib.de}. Supported by the DFG Research
Center \textsc{Matheon} "Mathematics for key technologies" in
Berlin.} \footnotetext[5]{Departamento de Matem\'atica Aplicada,
Universidad de Valladolid, Valladolid, Spain.
 ~E-mail: {\tt marial@mac.cie.uva.es}. Supported by DGI-MCYT under
 project MTM 2004-07194 cofinanced by FEDER funds.}
\footnotetext[7]{Mathematisches Institut,
  Universit\"at T\"ubingen,
  Auf der Morgenstelle 10,
  D--72076 T\"ubingen,
  Germany.
  ~
  E-mail:  {\tt lubich@na.uni-tuebingen.de}. Supported by DFG, SFB 382.}
\renewcommand{\thefootnote}{\arabic{footnote}}

\begin{abstract} We give an algorithm to compute $N$ steps of a convolution
quadrature approximation to a continuous temporal convolution
using only $O(N\, \log N)$ multiplications and  $O(\log N)$ active memory.
The method does not require evaluations of the convolution kernel,
but instead $O(\log N)$ evaluations of its Laplace transform,
which is assumed sectorial.
 The algorithm
can be used for the stable numerical solution with quasi-optimal
complexity of linear and nonlinear integral and
integro-differential equations of convolution type.
In a numerical example we apply it to solve a subdiffusion
equation with transparent boundary conditions.
\end{abstract}


\section{Introduction}
In this paper we give a fast and memory-saving algorithm for
computing the approximation of a continuous convolution
(possibly matrix $\times$ vector)
\begin{equation}\label{conv}
  \int_0^t f(t-\tau) \, g(\tau) \, d\tau~, \qquad 0\le t \le T,
\end{equation}
by a convolution quadrature with a step size $h>0$,
\begin{equation}\label{cq}
  \sum_{j=0}^n \omega_{n-j} \, g(jh)~, \qquad n=1,\dots,N,
\end{equation}
where the convolution quadrature weights $\omega_n$ are determined
from their generating power series as (see \cite{Lu88a, Lu88b, Lu04})
\begin{equation}\label{omega}
\sum_{n=0}^\infty \omega_n \zeta^n = F\Bigl( {\delta(\zeta)\over h} \Bigr).
\end{equation}
Here $F(s)$ is the Laplace transform of the (possibly matrix-valued)
convolution kernel $f(t)$, and $\delta(\zeta)=1-\zeta$
or $\delta(\zeta)=(1-\zeta)+\half(1-\zeta)^2$ for
the methods based on the first or second-order backward difference formula,
respectively. We will also consider a similar approximation based on
implicit Runge-Kutta formulas such as the Radau IIA
methods \cite{LuO93}.
Attractive features of such convolution quadratures are that they work well
for singular kernels $f(t)$, for
kernels with multiple time scales, and in situations where
only the Laplace transform $F(s)$ but not the
convolution kernel $f(t)$ is known analytically. Perhaps most importantly,
they enjoy excellent stability properties when used for the discretization of
integral equations or integro-differential equations of convolution type,
in a way often strikingly opposed to discretizations with more straightforward
quadrature formulas (see references in \cite{Lu04}).

The direct way to compute (\ref{cq}) is to first compute and store
the (possibly matrix-valued) weights
$\omega_0,\dots,\omega_N$, which can be done accurately with $O(N)$
evaluations of the Laplace transform $F(s)$ \cite{Lu88b}, and then
to compute the discrete convolution. Done naively, this requires
$O(N^2)$ multiplications (possibly matrix $\times$ vector)
and $O(N)$ active memory for the values $g(jh)$
and for the weights. Using FFT, the number of multiplications can be reduced to
$O(N\,\log N)$, and to $O(N\,(\log N)^2)$
in the case of integral equations where the values of $g(t)$ are not
known beforehand, but where $g(nh)$ is computed only in the $n$th time step
\cite{HaLS85}. However, that approach does not reduce the number of
$F$-evaluations and the memory requirements.

Here we give an algorithm, also applicable in the case of
linear and nonlinear integral
equations, which computes (\ref{cq}) in a way that requires
\begin{itemize}
\item $O(N\, \log N)$ multiplications,
\item $O(\log N)$ evaluations of the Laplace
transform $F(s)$, and
\item $O(\log N)$ active memory.
\end{itemize}
The history $g(jh)$ for $j=0,\dots,N$ is forgotten in this
algorithm, and only logarithmically
few linear combinations of the $g$-values are kept in memory.
These are
obtained by solving numerically, with step size $h$,
initial value problems of the
form $y'=\lambda y + g$ with complex $\lambda$.
The weights $\omega_n$ ($n\le N$) are not computed explicitly,
except the first few, e.g., the first 10  weights.

The algorithm presented here uses ideas of the fast convolution
algorithm of~\cite{LuS02}, which instead of (\ref{cq}) makes a
different approximation to the continuous convolution.
The stability properties of the second-order method of \cite{LuS02}
in integro-differential equations such as those of Section~5
are, however, extremely difficult to analyze (cf.~also \cite{Sch02})
and remain entirely unclear for higher-order extensions.
Here we show how the convolution quadratures (\ref{cq}) with all their known
favorable properties can be implemented in an equally fast and
memory-saving way.

Following the error analysis of \cite{LoLPS,LoP03} we give exponentially convergent error
bounds for the contour integral approximations that are
employed in this algorithm. They ensure that the constants hidden in
the $O$-symbols of the above work estimates depend only logarithmically
on the error tolerance for these contour integral approximations.

We assume a {\it sectorial} Laplace transform $F(s)$:
\begin{equation}\label{sector}
  \begin{array}{c}
    \hbox{$F(s)$ is analytic in a sector $|\arg(s-c)|<\pi-\varphi$ with
      $\varphi<\half\pi$, % and real $c$,
      and there}
    \\[2mm]
    |F(s)| \le M\, |s|^{-\nu}\quad\hbox{for some real $M$ and $\nu>0$}.
  \end{array}
\end{equation}
The inverse Laplace transform is then given by
\begin{equation}\label{inv}
  f(t) = {1\over 2\pi i} \int_\Gamma   e^{t\lambda }\,F(\lambda)\,  d\lambda,
  \qquad t>0,
\end{equation}
with $\Gamma$ a contour in the sector of analyticity, going to infinity
with an acute angle to the negative real half-axis and
oriented with increasing imaginary part. The function $f(t)$
is analytic in $t>0$ and satisfies
\begin{equation}\label{fbound}
  |f(t)|\le C \,t^{\nu-1}\,e^{ct},\qquad t>0,
\end{equation}
and is therefore locally integrable.
(The absolute values on the left-hand sides of the bounds (\ref{sector}) and
(\ref{fbound}) are to be interpreted as matrix norms for
matrix-valued convolution kernels.)

In Section~2 we review convolution quadrature based on
multistep and Runge-Kutta methods. We give a contour integral representation
of the convolution quadrature weights whose discretization along
hyperbolas or Talbot contours is discussed in Section~3.
The fast and oblivious convolution algorithm is formulated in
Section~4. Finally, in Section~5 we give the results of numerical
experiments with integral and integro-differential equations
originating from regular and anomalous diffusion problems.

\pagebreak[3]

\section{Convolution quadrature}
In this section we review briefly convolution quadrature
%based on linear multistep and Runge-Kutta methods
and give a
contour integral representation of the convolution quadrature
weights on which the fast algorithm of this paper is based.

\subsection{Convolution quadrature based on multistep methods}
\label{subsec:ms}
We consider the convolution quadrature (\ref{cq}) with weights (\ref{omega}).
By (\ref{sector}) and Cauchy's integral formula we have, with a contour
$\Gamma$ as in (\ref{inv}),
$$
\sum_{n=0}^\infty \omega_n \zeta^n = F\Bigl( {\delta(\zeta)\over h} \Bigr)=
{1\over 2\pi i}
\int_{\Gamma}
\Bigl( {\delta(\zeta)\over h} -\lambda\Bigr)^{-1} \,F(\lambda) \, d\lambda.
$$
Hence, with $e_n(z)$ defined by
\begin{equation}\label{en-def}
  (\delta(\zeta)-z)^{-1} = \sum_{n=0}^\infty e_n(z) \,\zeta^n ,
\end{equation}
we have the integral formula
\begin{equation}\label{om-int}
  \omega_n = {h\over 2\pi i}
  \int_\Gamma    e_n(h\lambda)\, F(\lambda) \, d\lambda,
\end{equation}
which can be viewed as the discrete analog of (\ref{inv}).
For the backward Euler discretization $\delta(\zeta)=1-\zeta$ we note
the explicit formula
\begin{equation}\label{en-bdf1}
  e_n(z) = (1-z)^{-n-1},
\end{equation}
which is of the form
$e_n(z) = q(z) r(z)^n$ with  %\qquad\hbox{with}\quad
$r(z) = \frac1{1-z}$ and $  q(z) =\frac1{1-z}$.

For the second-order BDF method, where
$\delta(\zeta)=\sum_{k=1}^p \frac1k (1-\zeta)^k$ with $p=2$,
we obtain from a partial fraction
decomposition of $(\delta(\zeta)-z)^{-1}$ that
\begin{equation}\label{en-bdf2}
  e_{n}(z) = \frac{1}{\sqrt{1+2z}}
  \Bigl((2-\sqrt{1+2z})^{-n-1} - (2+\sqrt{1+2z})^{-n-1}\Bigr),
\end{equation}
which is of the form $e_n(z) = q_1(z) r_1(z)^n + q_2(z) r_2(z)^n$.
Connoisseurs of Cardano's formulas find analogous formulas to (\ref{en-bdf2})
also for
the BDF methods of orders 3 and~4.

\subsection{Convolution quadrature based on Runge-Kutta methods}
\label{subsec:rk}
We consider an implicit Runge-Kutta method with
coefficients $a_{ij}$, $b_j$, $c_i$
for $i,j=1,\dots,m$. We denote the Runge-Kutta matrix by
$\A=(a_{ij})$, the row vector of the weights by $b^T=(b_j)$, and
the stability function by
$$
r(z) = 1 + z b^T(I-z\A)^{-1}\1,
$$
where $\1=(1,\dots,1)^T$. We assume that all
eigenvalues of the Runge-Kutta matrix $\A$ have positive real part
and, for simplicity, that the method is A-stable
and the row vector of the weights
equals the last line of the Runge-Kutta matrix,
$$
b_j = a_{mj} \quad\hbox{ for }\ j=1,\dots,m,
$$
and correspondingly $c_m=1$.
These conditions are in particular satisfied by the Radau IIA
family of Runge-Kutta methods \cite{HaW96}.
From such a Runge-Kutta method, a convolution quadrature is constructed
as follows \cite{LuO93}: Let
\begin{equation}\label{Delta}
  \Delta(\z)=\Bigl(\A + {\z \over 1-\z}\1 b^T\Bigr)^{-1}
\end{equation}
and define weight matrices $W_n$ by
\begin{equation}\label{Wn}
  \sum_{n=0}^\infty W_n\z^n = F\Bigl( {\Delta(\z) \over h}\Bigr).
\end{equation}
Let $\omega_n=(\omega_n^1,\dots,\omega_n^m)$ denote the last row of $W_n$.
Then an approximation to the convolution integral (\ref{conv}) at time
$t_{n+1}=(n+1)h$
is given by
\begin{equation}\label{rk-cq}
 u_{n+1}=   \sum_{j=0}^n \sum_{i=1}^m  \omega_{n-j}^i\, g( t_j + c_ih) =
   \sum_{j=0}^n  \omega_{n-j}\, g_j
\end{equation}
with the column vector $g_j = \bigl( g(t_j+c_ih) \bigr)_{i=1}^m$.
For a Runge-Kutta method of classical order $p$ and stage order $q$, this
approximation is known to be
convergent of the order $\min(p,q+1+\nu)$ with $\nu$ of
(\ref{sector}).

With the row vector $e_n(z)=(e_n^1(z),\dots,e_n^m(z))$
defined as the last row of
the $m\times m$ matrix $E_n(z)$ given by
\begin{equation}\label{En-def}
  (\Delta(\zeta)-zI_m)^{-1} = \sum_{n=0}^\infty E_n(z) \,\zeta^n ,
\end{equation}
we obtain an integral formula like in (\ref{om-int}),
\begin{equation}\label{W-int}
  \omega_n = {h\over 2\pi i}
  \int_\Gamma    e_n(h\lambda)\otimes F(\lambda)\, d\lambda.
\end{equation}
For $n\ge 0$, $\, e_n(z)$ is given as
\begin{equation}\label{en-rk}
e_n(z) = r(z)^n q(z)
\end{equation}
with the row vector $q(z) =  b^T (I-z\A)^{-1}$;
cf.~Lemma~2.4 in~\cite{LuO93}.
We note that
\begin{equation}\label{yrk}
y_{n+1}^{(\lambda)} = h\sum_{j=0}^n e_{n-j}(h\lambda)\, g_j
\end{equation}
is the Runge-Kutta approximation at time $t_{n+1}$ of the linear
initial value problem
\begin{equation}\label{y-lambda}
y'=\lambda y + g(t), \quad y(0)=0\,.
\end{equation}
The convolution quadrature (\ref{rk-cq})
is thus interpreted as
$$
 u_{n+1} = {1\over 2\pi i}
  \int_\Gamma   F(\lambda) \, y_{n+1}^{(\lambda)} \, d\lambda\; ;
$$
see Proposition~2.1 in \cite{LuO93}.

\section{Approximation of the contour integrals}
\label{Sec.Contour}
The fast convolution algorithm will be based on discretizing
the integrals in (\ref{om-int}) and (\ref{W-int})
along suitable complex
contours. This approximation is discussed in the present section.

\subsection{Quadrature on Talbot contours and hyperbolas}
The fast algorithm  approximates
the quadrature weights $\omega_n$
by linear combinations of the exponential approximations $e_n(h\lambda)$,
locally on a sequence of fast-growing
time intervals $nh \in I_\ell$:
\begin{equation}\label{intervals}
  I_\ell = [B^{\ell-1}h, 2B^{\ell}h),
\end{equation}
where the base $B>1$ is an integer. For example,
$B=10$ was found a good choice in our numerical experiments.
The approximation  on $I_\ell$ results
from applying the trapezoidal rule to
a parametrization of the contour integral for the convolution quadrature
weights,
\begin{equation}\label{eq:num-int}
  \omega_n = {h\over 2\pi i} \int_{\Gamma_\ell}
  e_n(h\lambda)\otimes F(\lambda) \, d\lambda
  \approx
 h\sum_{k=-K}^{K} w_k^{(\ell)}
 e_n(h\lambda_k^{(\ell)})\otimes F(\lambda_k^{(\ell)})\, ,
  \quad\ nh\in I_\ell,
\end{equation}
with an appropriately chosen complex contour $\Gamma_\ell$.
The number of quadrature points on $\Gamma_\ell$, $2K+1$,
is chosen independent
of $\ell$. It is much smaller than what would be required for a uniform
approximation of the contour integral on the whole interval $[0,T]$.
Only a few of the first convolution quadrature weights, $\omega_n$
for $n\le N_0$ (e.g., $N_0=10$), are
approximated differently, using the trapezoidal rule discretization of the
integral over a circle as discussed in \cite{Lu88b,LuO93}:
\begin{equation}\label{w-circle-int}
\omega_n = \hbox{ last row of }\, {1\over 2\pi i} \int_{|\zeta|=\rho}
\z^{-n-1}\, F\Bigl( {\Delta(\z)\over h}\Bigr)\, d\z,
\qquad n\le N_0.
\end{equation}
The numerical integration in (\ref{om-int}) or (\ref{W-int})
is done by applying the
trapezoidal rule with equidistant steps to a parameterization of a
hyperbola~\cite{LoP03} or a Talbot contour~\cite{Talbot79,Rizzardi95}.
\begin{figure}[ht]
 \psfrag{lambda }[cc][bl]{$\mu$}
 \psfrag{sigma }[cc][cc]{$\sigma$}
 \psfrag{lsnh }[cc][cB]{$\frac{\mu\nu\pi}{2}$}
 \psfrag{lsn }[cc][cr]{$\mu \nu \pi$}
 \includegraphics[width= 0.49\textwidth, height=0.2\textheight]{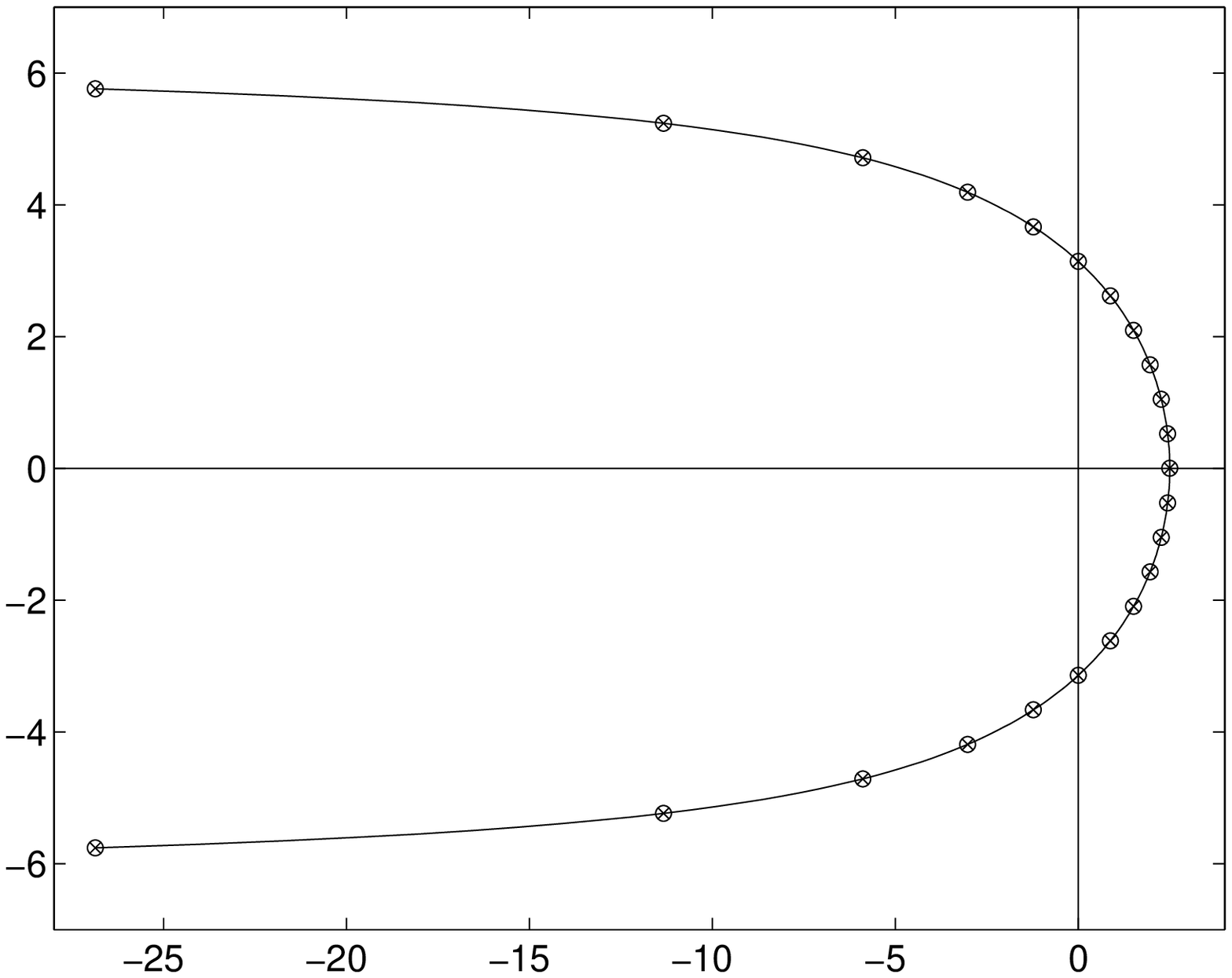}
 \includegraphics[width= 0.49\textwidth, height=0.2\textheight]{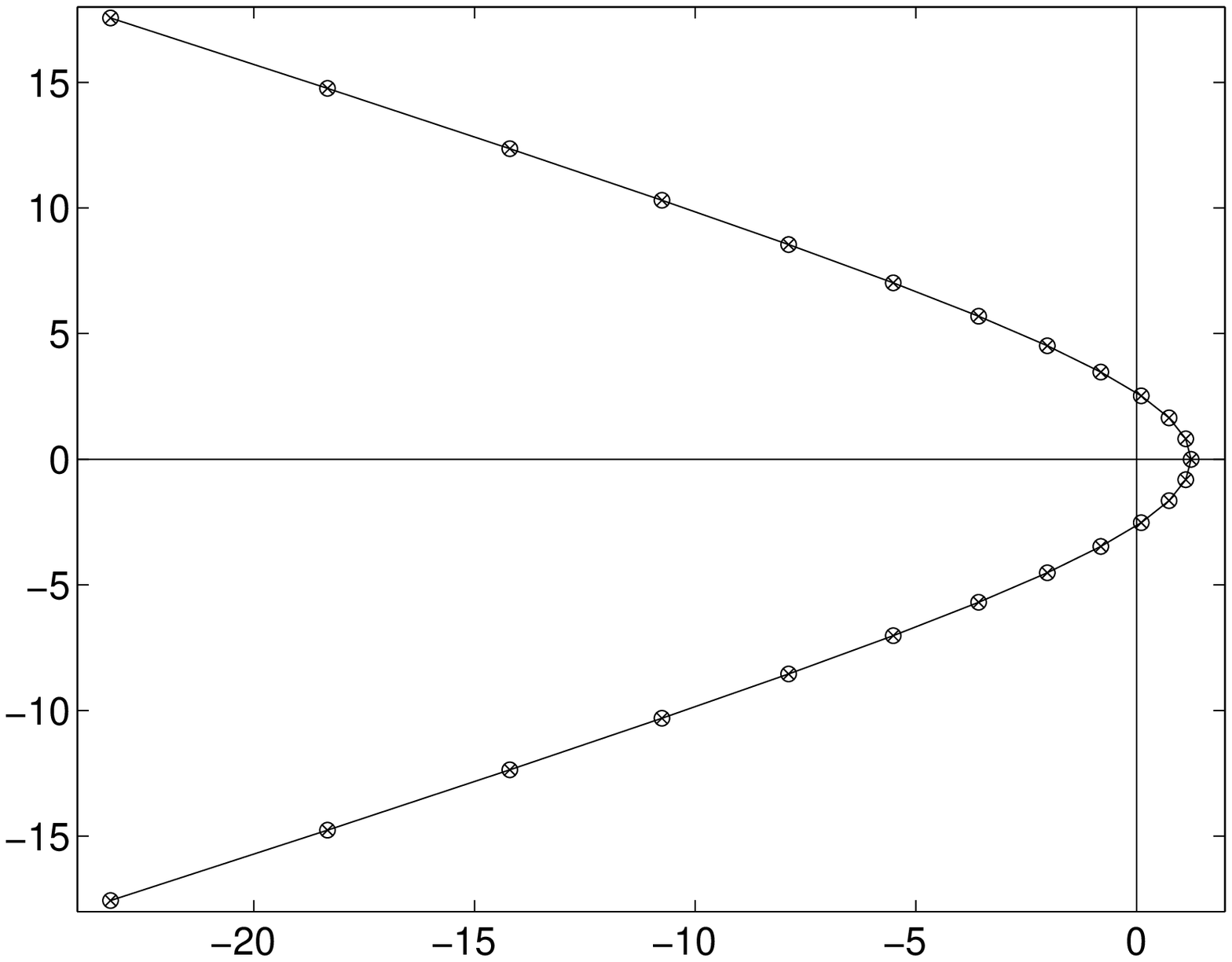}
 \caption{Talbot contour (left) and hyperbola (right).}
\label{Fig.talbotcont}
\end{figure}
The Talbot contour is given by
\begin{eqnarray}
  (-\pi,\pi) \to \Gamma &:& \theta \mapsto \gamma(\theta) = \sigma +
  \mu \left(
    \theta \cot(\theta) + i \nu \theta \right)
  \label{Eq.contourparametrisierung}
\end{eqnarray}
where the parameters $\mu$, $\nu$ and $\sigma$ are such that the
singularities of $F(s)$ lie to the left of the contour
and that the singularities of $e_n(hs)$ lie to the right of the contour.
See left part of Figure~\ref{Fig.talbotcont} for $\sigma = 0$.
The parameter $\mu$ will depend on $\ell$ via the right end-point of $I_\ell$,
which yields a Talbot
contour $\Gamma_\ell$ depending on the approximation interval~$I_\ell$.
The weights and quadrature points in
(\ref{eq:num-int}) are given by (omitting $\ell$ in the notation)
$$
w_k = -{i\over 2(K+1)}\: \gamma'(\theta_k)~, \quad\
\lambda_k = \gamma(\theta_k) \quad\ \mbox{ with }\quad
\theta_k={k\pi\over K+1}~.
$$
Alternatively, the hyperbola  is given by
\begin{eqnarray}
  \Rset \to \Gamma &:& \theta \mapsto \gamma(\theta) =
  \mu (1 - \sin(\alpha+i\theta))
  \label{Eq.contourparametrisierunghyperbola}
\end{eqnarray}
where the parameters $\mu$ and $\alpha$ are such that the
singularities of $F(s)$ lie to the left of the contour.
See the right part of Figure~\ref{Fig.talbotcont} for $\alpha = \pi/2 - 1/2$.
The weights and quadrature points in (\ref{eq:num-int}) are given by
(omitting $\ell$ in the notation)
$$
w_k = {i \tau \over 2 \pi}\: \gamma'(\theta_k)~, \quad
\lambda_k = \gamma(\theta_k) \quad\ \mbox{ with }\quad
\theta_k = k\tau~,
$$
where $\tau$ is a step length parameter.

\subsection{Numerical experiments}\label{subsec:numexp}
In view of the  examples of
Section~\ref{NumEx} we present here numerical experiments with
$$
f(t) = {1\over \sqrt{\pi t}}, \quad\ \hbox{for which }\quad
F(s)=s^{-1/2}.
$$

\begin{figure}[!h]
\centerline{\psfig{file=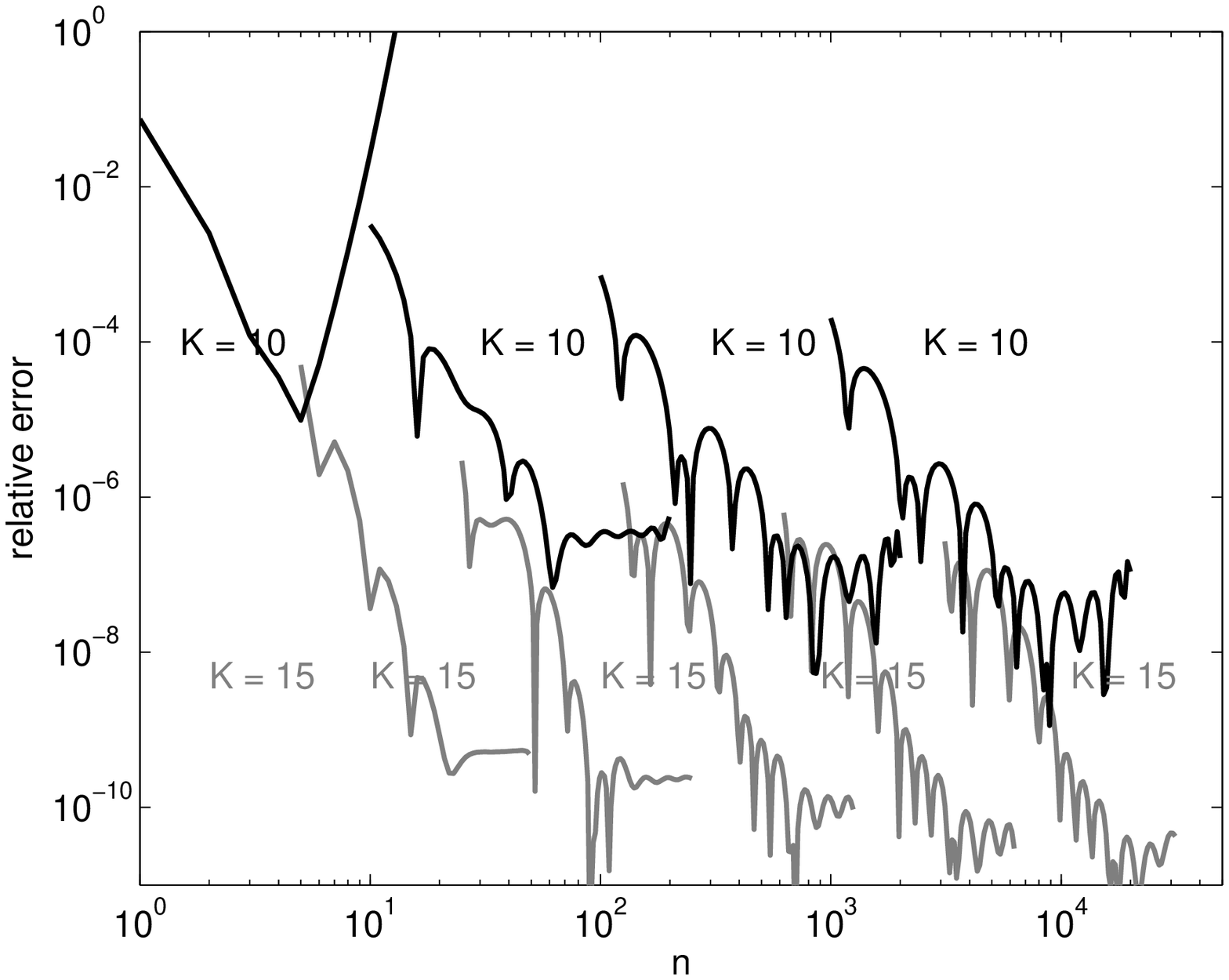,width=0.40\textwidth}
            \psfig{file=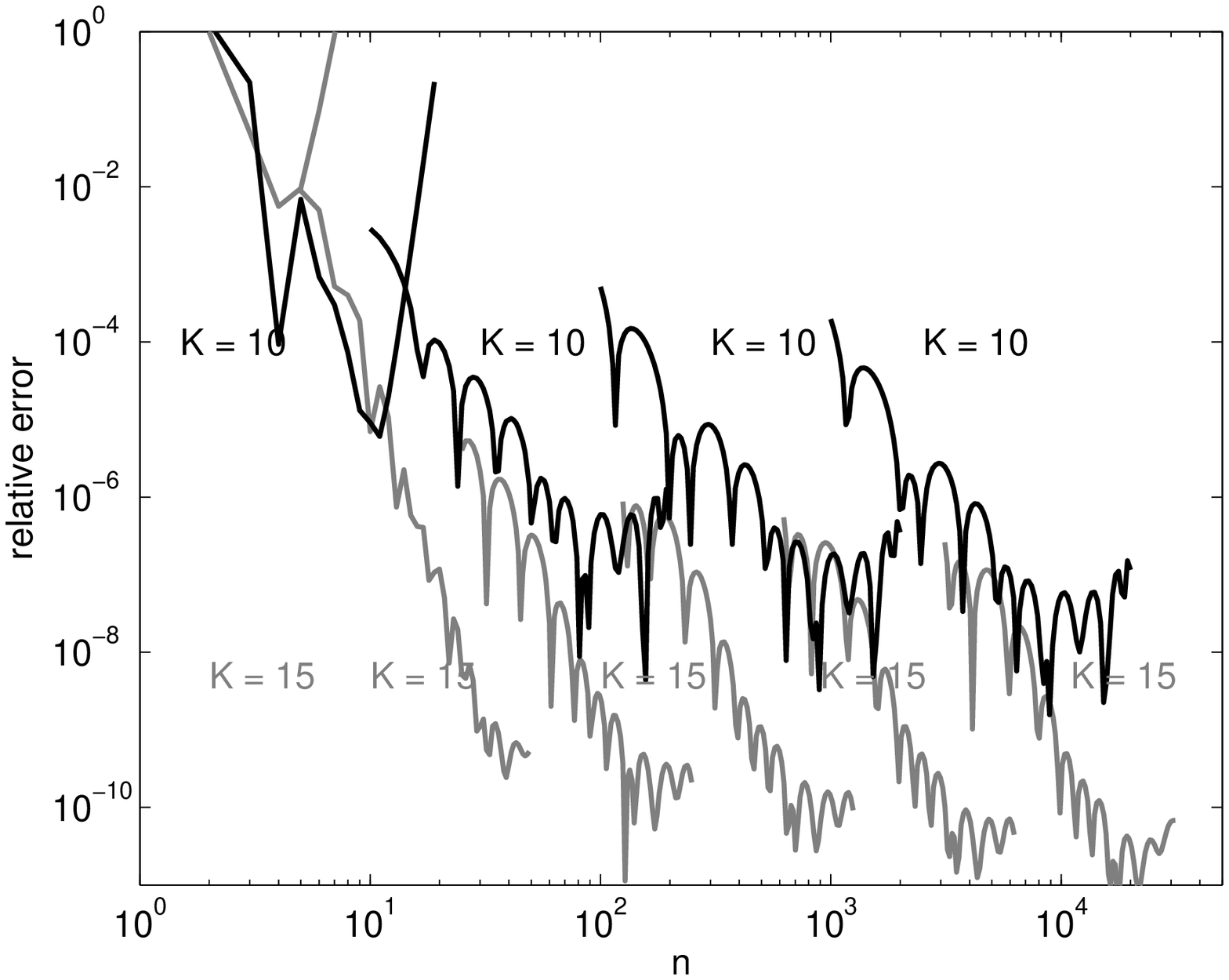,width=0.40\textwidth}}
\centerline{\psfig{file=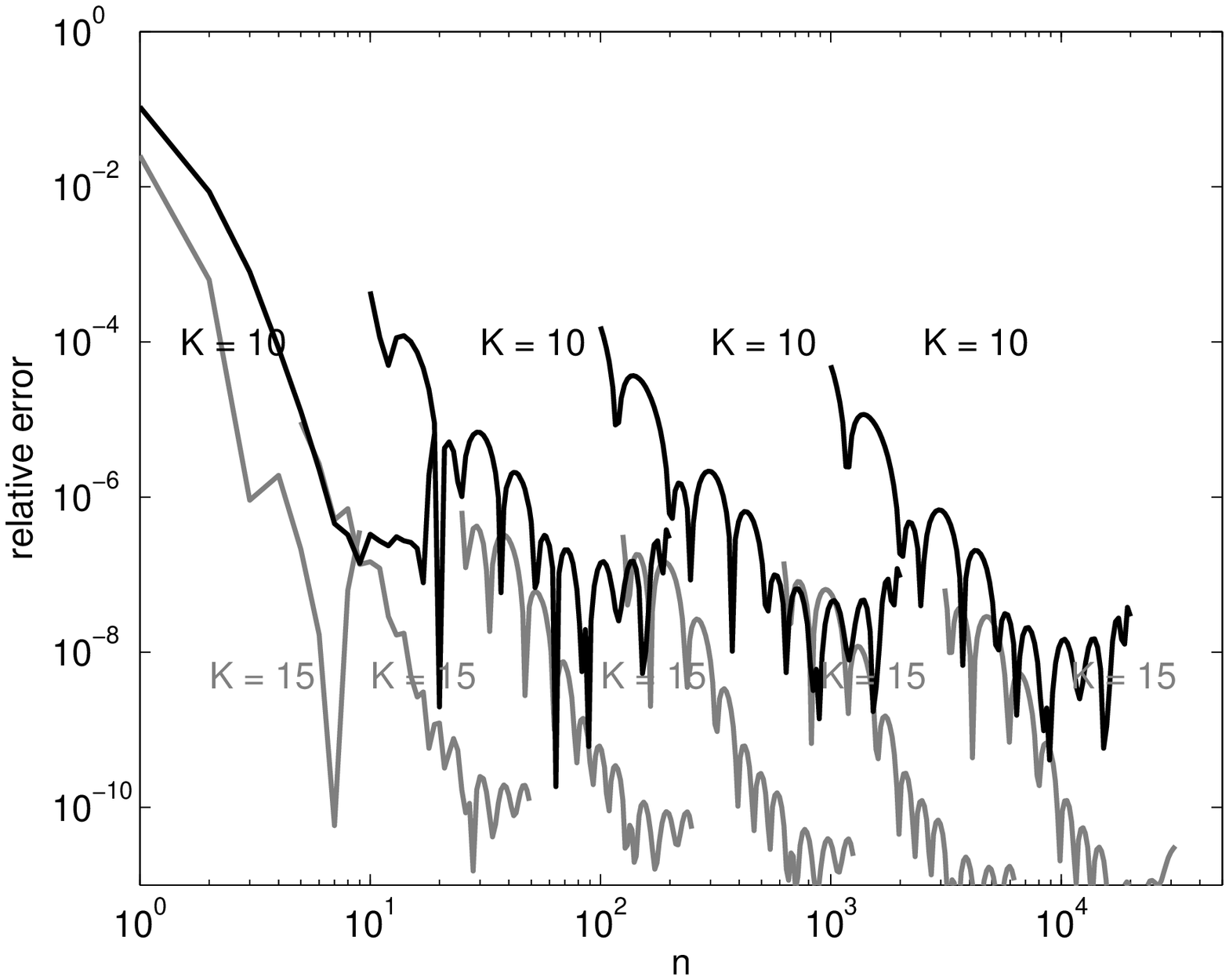,width=0.40\textwidth}
            \psfig{file=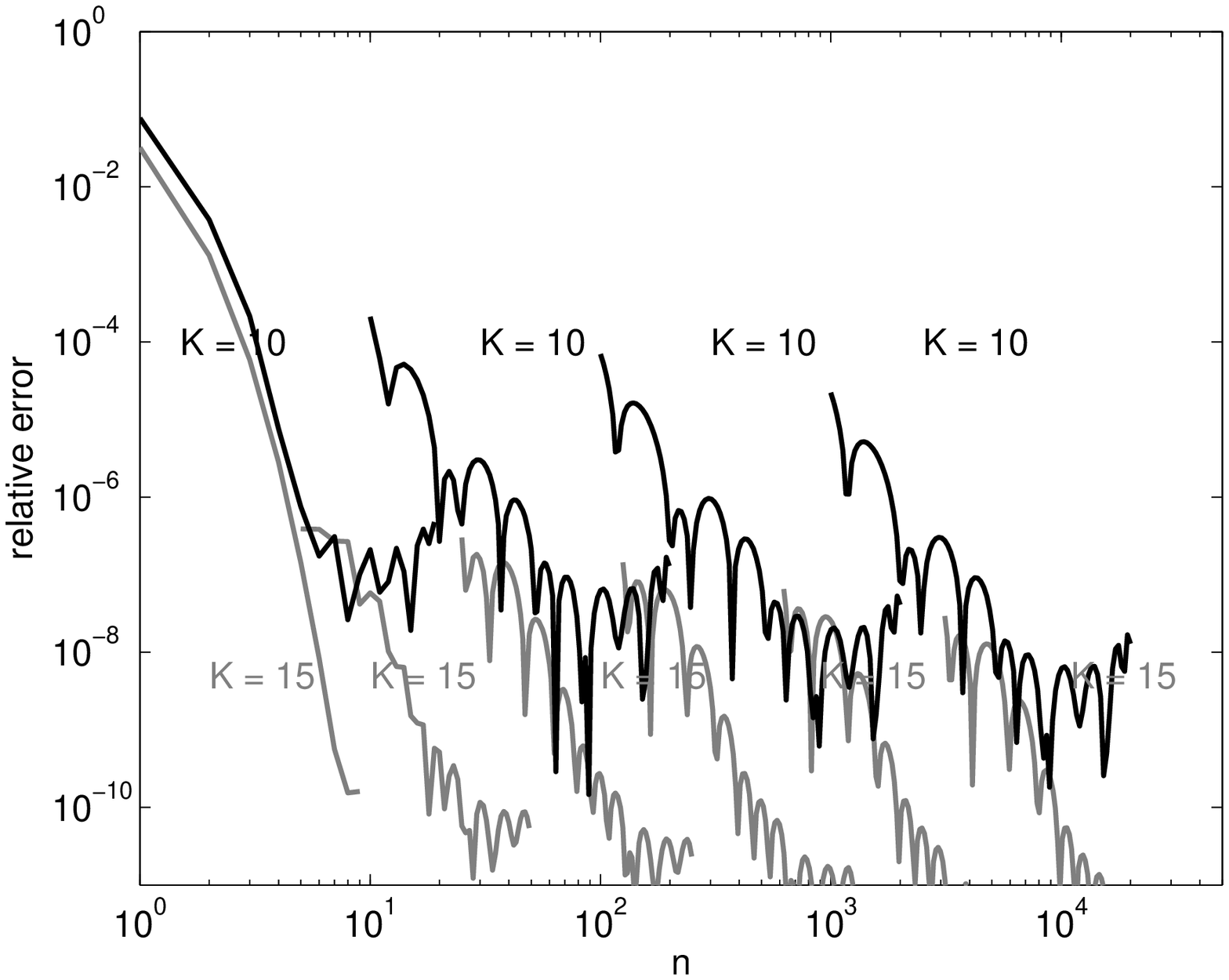,width=0.40\textwidth}}
\caption{Talbot quadrature errors versus time for $K = 15$, $B=5$ and $K= 10$, $B=10$
  for different Integrators. (Implicit Euler, BDF(2), RadauIIA(3) and RadauIIA(5)
  in clockwise order starting from the upper left corner)}
\label{fig:talbot-error}
\end{figure}

\begin{figure}[!h]
\centerline{\psfig{file=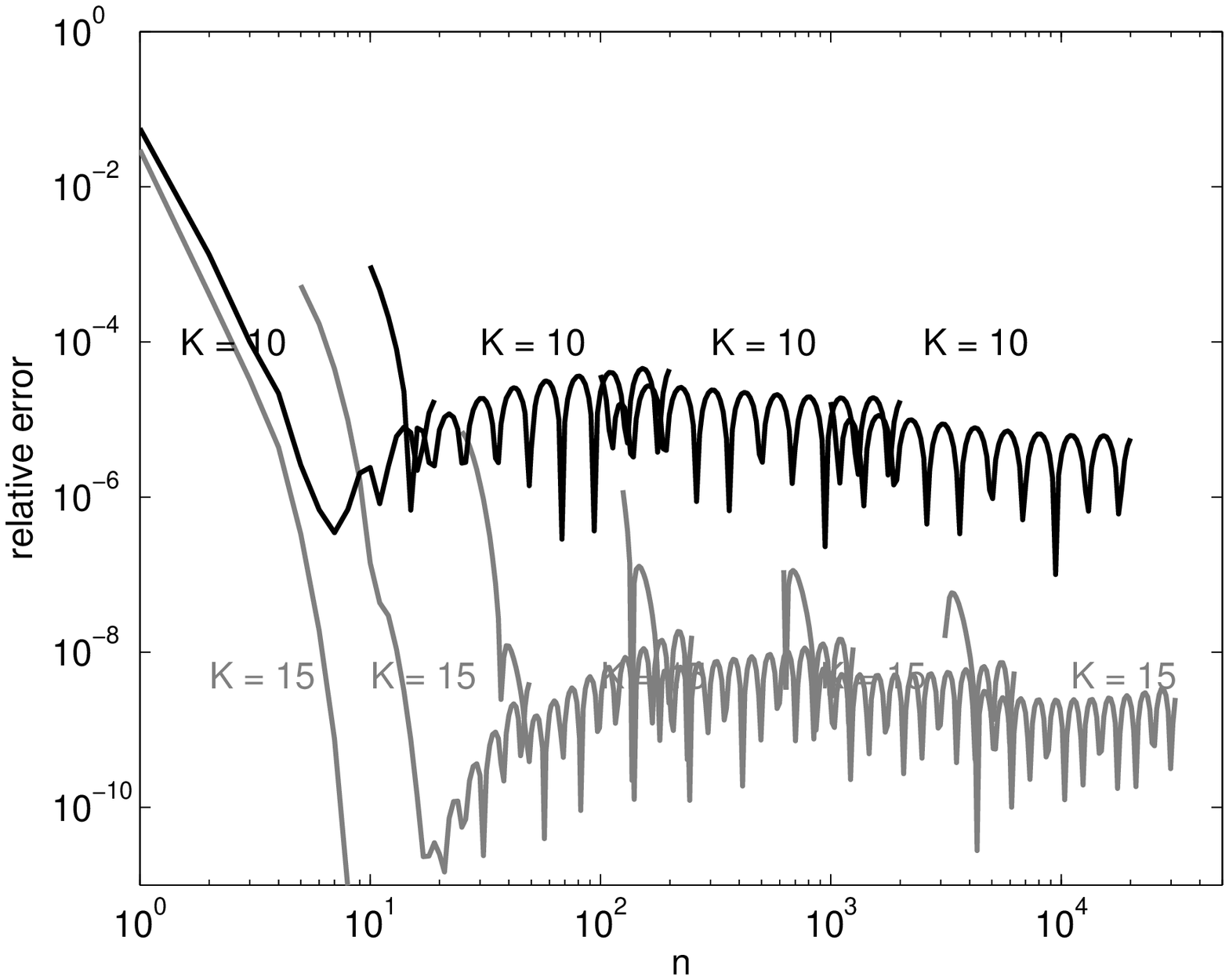,width=0.40\textwidth}
            \psfig{file=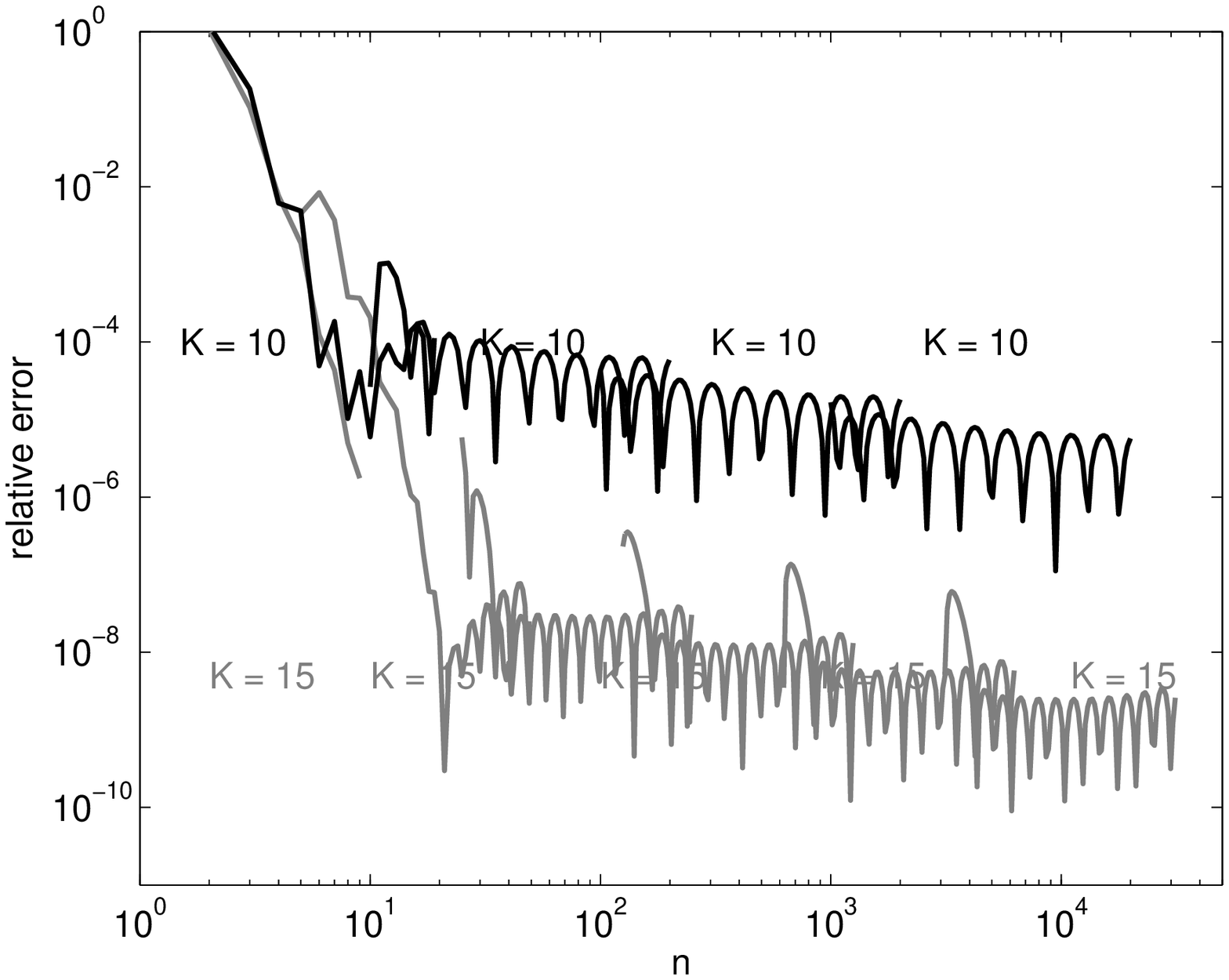,width=0.40\textwidth}}
\centerline{\psfig{file=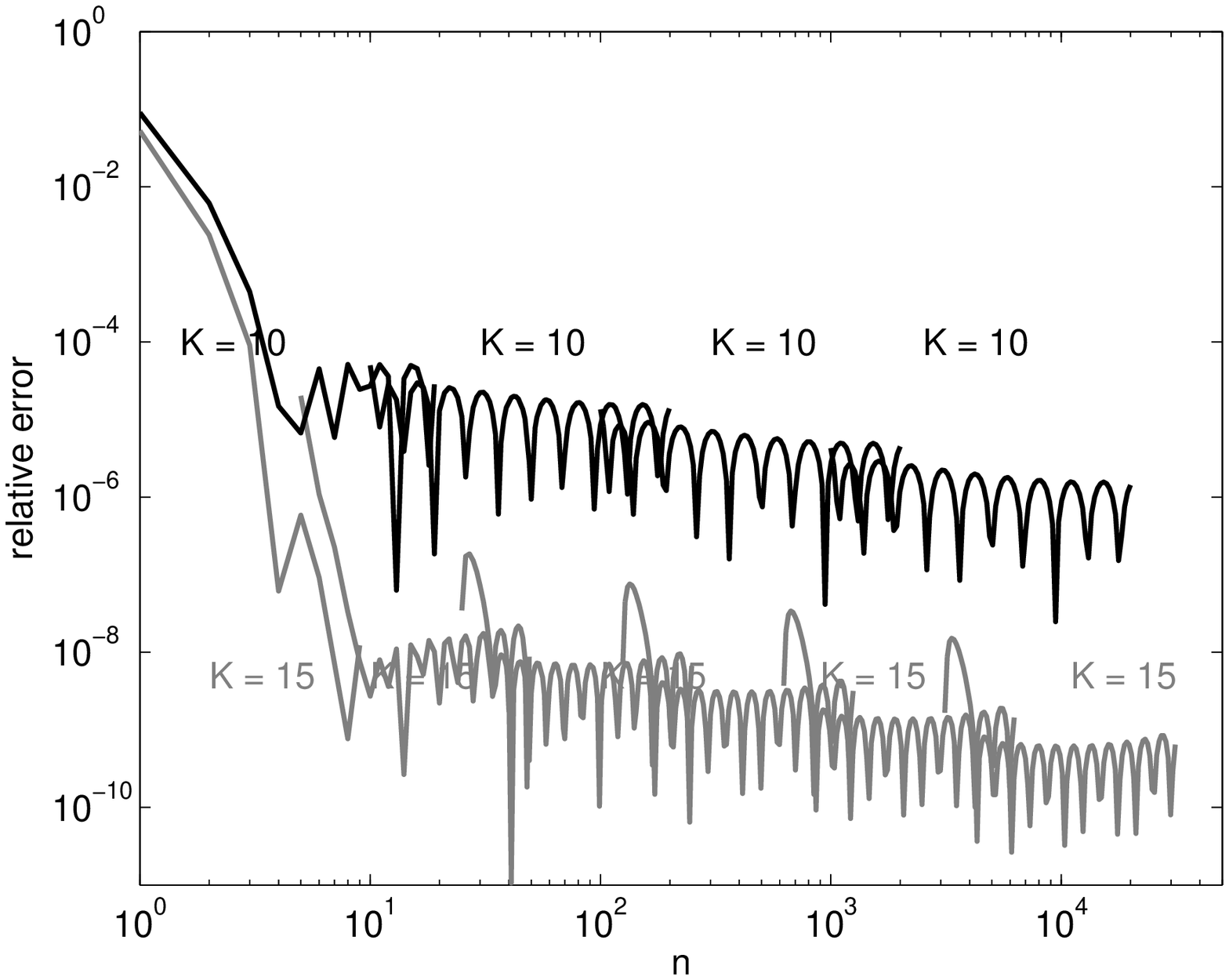,width=0.40\textwidth}
            \psfig{file=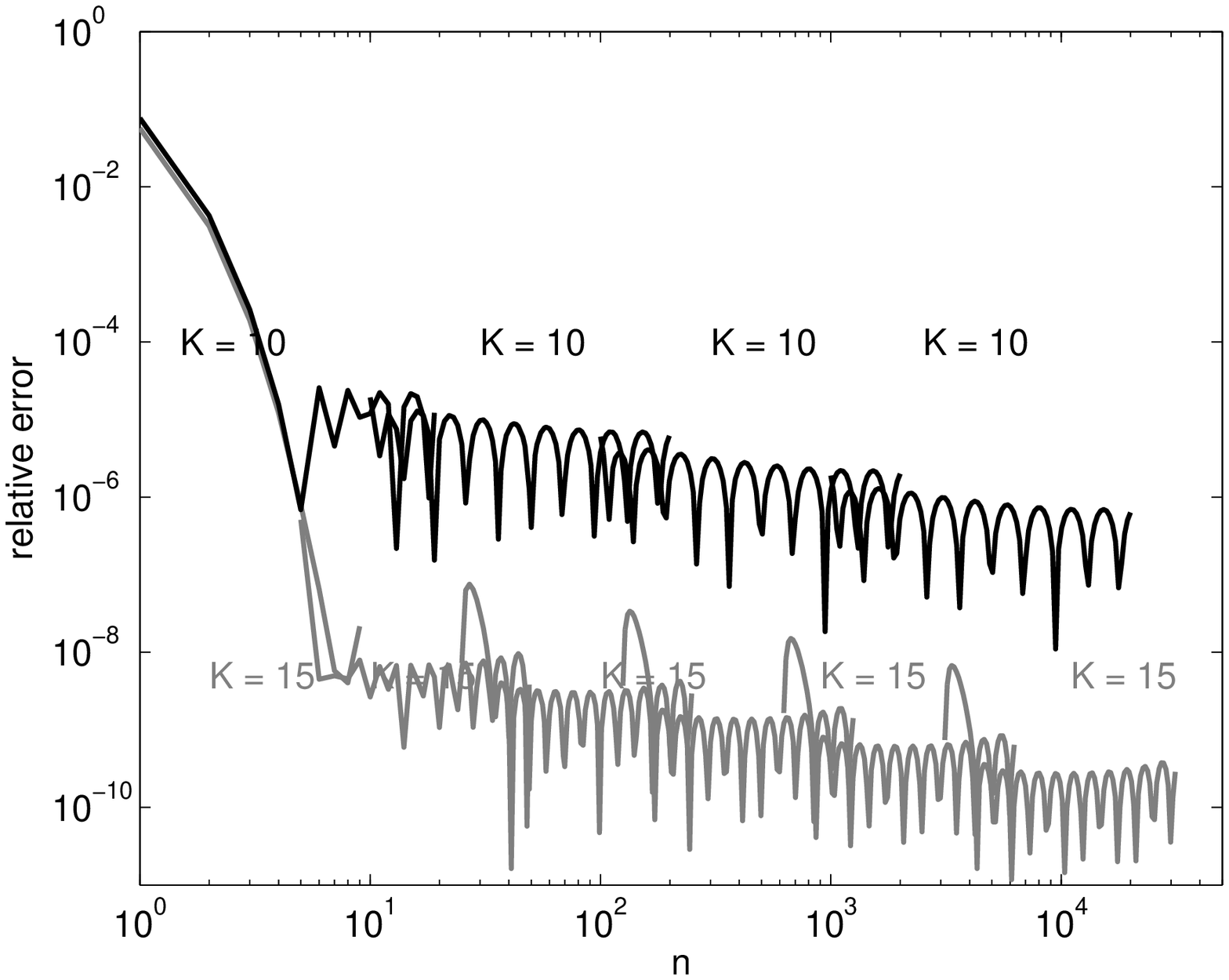,width=0.40\textwidth}}
\caption{Hyperbola quadrature errors versus time for $K = 15$, $B=5$ and $K= 10$, $B=10$
  for different Integrators. (Implicit Euler, BDF(2), RadauIIA(3) and RadauIIA(5)
  in clockwise order starting from the upper left corner)}
\label{fig:hyperbola-error}
\end{figure}

%\textit{For the $\rho/(s(\rho + s))$ example Talbot error and Hyperbola error are comparable.}

%%   \begin{figure}[!h]
%%     \centerline{\psfig{file=TalboterrorN10N15_dis_RKVRadauIIA1Rover_Sx_RpS__.eps,width=0.40\textwidth}
%%       \psfig{file=TalboterrorN10N15_dis_BDF2Rover_Sx_RpS__.eps,width=0.40\textwidth}}
%%     \centerline{\psfig{file=TalboterrorN10N15_dis_RKVRadauIIA3Rover_Sx_RpS__.eps,width=0.40\textwidth}
%%       \psfig{file=TalboterrorN10N15_dis_RKVRadauIIA5Rover_Sx_RpS__.eps,width=0.40\textwidth}}
%%     \caption{Talbot quadrature errors versus time for $K = 15$, $B=5$ and $K= 10$, $B=10$
%%       for different Integrators. (Implicit Euler, BDF(2), RadauIIA(3) and RadauIIA(5)
%%       in clockwise order starting from the upper left corner)}
%%     \label{fig:talbot-error2}
%%   \end{figure}

%%   \begin{figure}[!h]
%%     \centerline{\psfig{file=HyperbolaerrorN10N15_dis_RKVRadauIIA1Rover_Sx_RpS__.eps,width=0.40\textwidth}
%%       \psfig{file=HyperbolaerrorN10N15_dis_BDF2Rover_Sx_RpS__.eps,width=0.40\textwidth}}
%%     \centerline{\psfig{file=HyperbolaerrorN10N15_dis_RKVRadauIIA3Rover_Sx_RpS__.eps,width=0.40\textwidth}
%%       \psfig{file=HyperbolaerrorN10N15_dis_RKVRadauIIA5Rover_Sx_RpS__.eps,width=0.40\textwidth}}
%%     \caption{Hyperbola quadrature errors versus time for $K = 15$, $B=5$ and $K= 10$, $B=10$
%%       for different Integrators. (Implicit Euler, BDF(2), RadauIIA(3) and RadauIIA(5)
%%       in clockwise order starting from the upper left corner)}
%%     \label{fig:hyperbola-error2}
%%   \end{figure}

The error is calculated with respect to a reference solution, obtained for
a discretization of the contour integral with a large number of
integration points.
For the Radau IIA methods of order~3 and~5, where the $\omega_{n}$ are
row vectors of dimension 2 and 3, respectively,
we plot the error of the last entry.

Using the Tabot contours, the following choices of parameters were
found to give good results.
A relative accuracy of about $ 10^{-3}$ on the interval $I_\ell$
for $\ell \ge 2$ with right end-point $T_\ell$ is obtained with
$B=10$, $K=10$, $\mu=8/T_\ell$, $\nu=0.6$.
For a relative approximation error of~$10^{-6}$, take $B=5$, $K=15$,
and the other parameters as before, cf.~Fig.~\ref{fig:talbot-error}.
For $n>20$ there is no substantial difference between the different
Runge-Kutta methods.
Since the approximations to the first few convolution quadrature weights
are poor, they will not be used in the algorithm.

Using the hyperbola contours, a relative accuracy of about
$ 10^{-4}$ on the interval $I_{\ell}$ for $\ell \ge 2$
with right end-point $T_{\ell}$ is obtained with
$B=10$, $K=10$, $\alpha = 1 $,
$\mu =3.6/T_\ell$ and $\tau = 0.64$.
%and $\mu$ and $\tau$ determined as above.
For a relative approximation error of~$3\cdot 10^{-8}$, we
take $B=5$, $K=15$, $\alpha = 1$,
cf.~Fig.~\ref{fig:hyperbola-error}.
For $n>20$ there is again no essential difference
between the different Runge-Kutta methods.

\begin{figure}[htbp]
%%   \centerline{\psfig{file=TalboterrorN10_N640_dis_RKVRadauIIA1oneoversqrt.eps,width=0.40\textwidth}
%%     \psfig{file=TalboterrorN10_N640_dis_BDF2oneoversqrt.eps,width=0.40\textwidth}}
%%   \centerline{\psfig{file=TalboterrorN10_N640_dis_RKVRadauIIA3oneoversqrt.eps,width=0.40\textwidth}
%%     \psfig{file=TalboterrorN10_N640_dis_RKVRadauIIA5oneoversqrt.eps,width=0.40\textwidth}}
\includegraphics[width=0.49\textwidth]{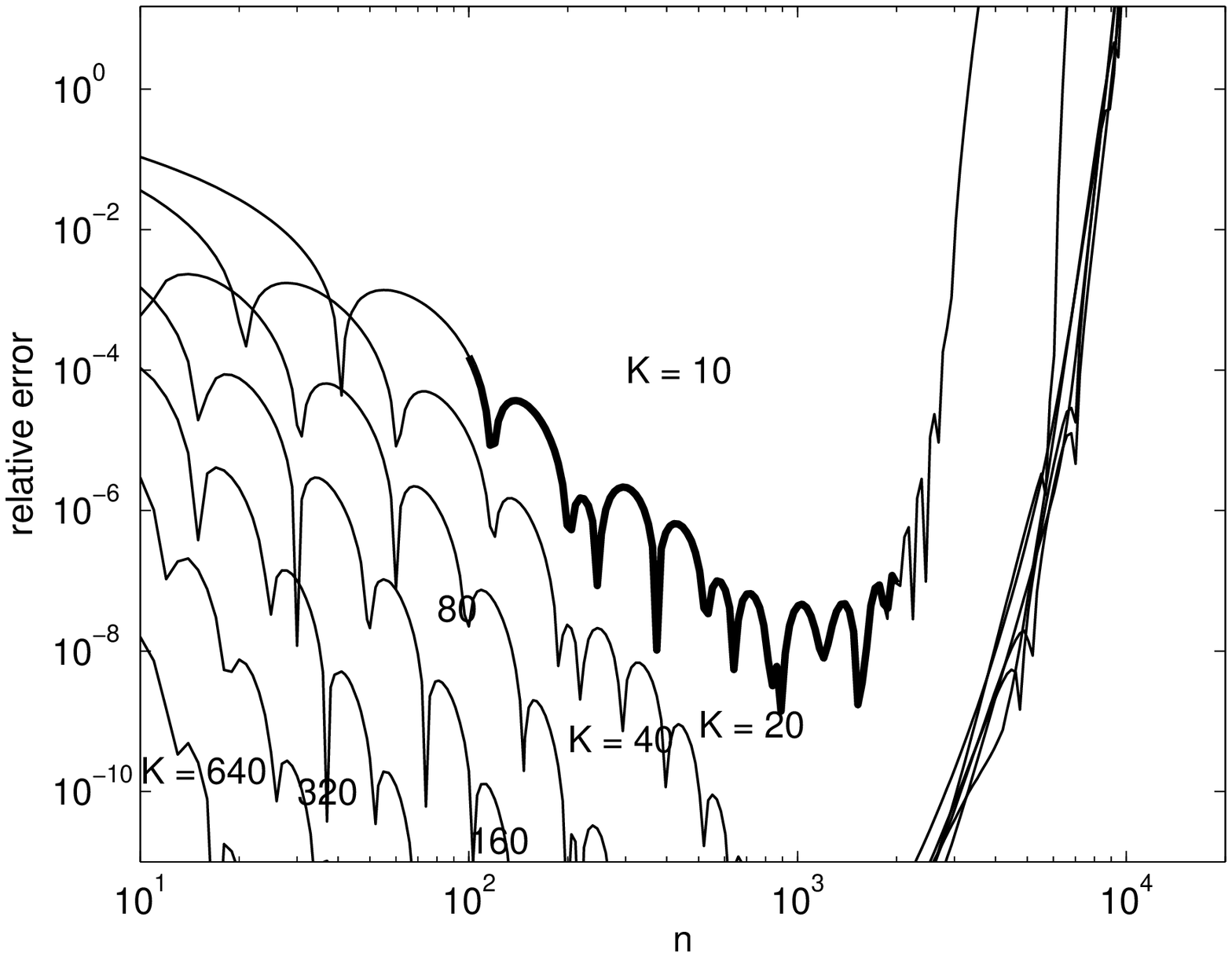}
\includegraphics[width=0.49\textwidth]{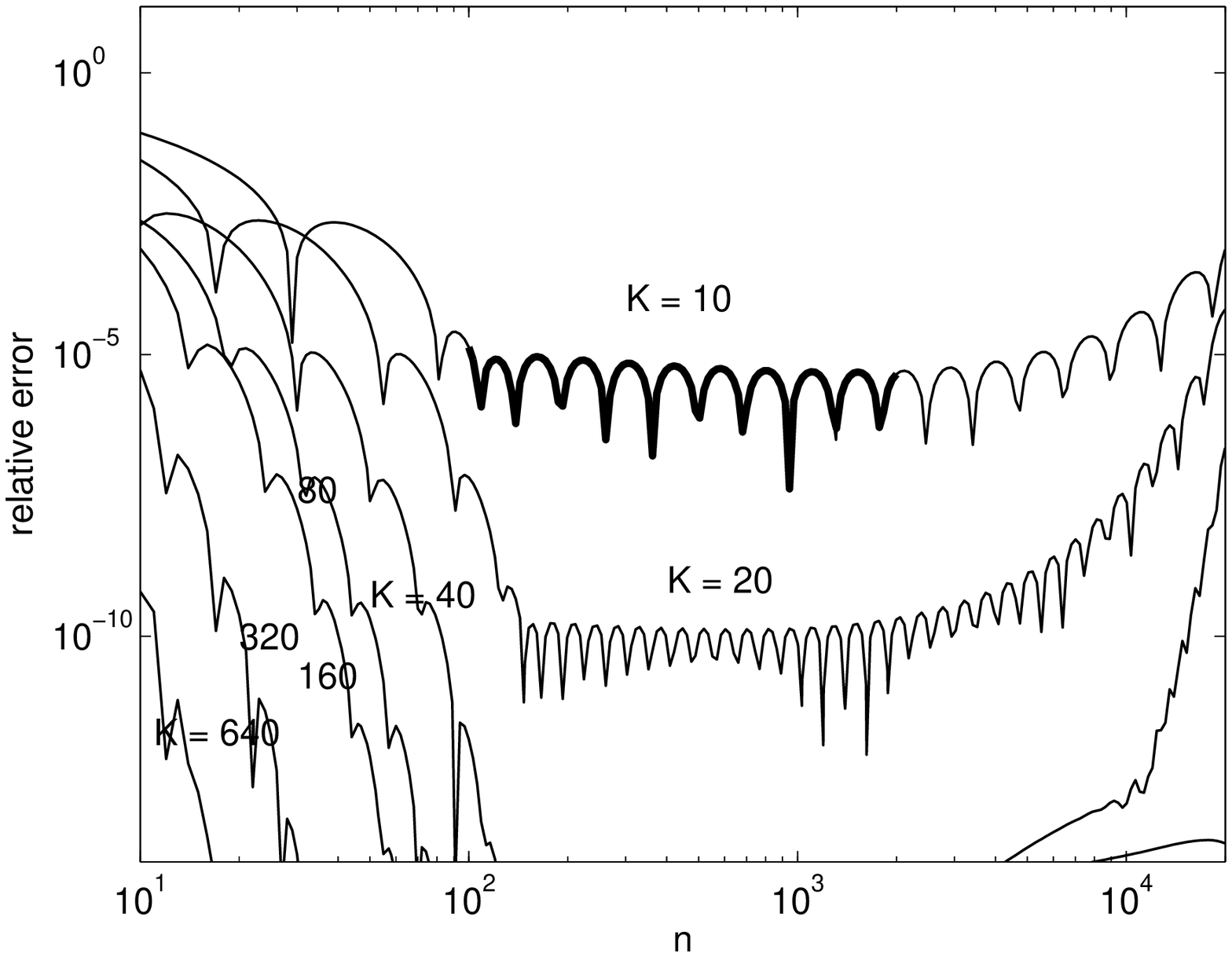}
% \centerline{\psfig{file=TalboterrorN10_N640_dis_RKVRadauIIA3oneoversqrt.eps,width=0.60\textwidth}}
  \caption{Talbot (left) and hyperbola (right) quadrature error versus $n$
  for the RadauIIA(3) method with
    $K = 10,20,40,80,160,320,640$.
    The bold parts of the error curve  correspond to
     the lower left parts of Figs.~\ref{fig:talbot-error} and
     \ref{fig:hyperbola-error}. Note the different scaling of the $y$ axis.}
  \label{fig:talbot-errora}
\end{figure}

Fig.~\ref{fig:talbot-errora} % and~\ref{fig:hyperbola-errora}
shows
the relative errors  on the interval
$[10,20000]$ (similar for any interval $[a,2000a]$ with $a>10$) for the
RadauIIA(3) method with $K=10,20,40,80,$ $160,320,640$. For the implicit Euler,
the BDF(2) and the RadauIIA(5) method these error plots look similar.
This behavior of the errors clearly
demonstrates the advantage of using local approximations.
With $B=10$, we need three approximation intervals
to cover the interval $[10,20000]$, so that for a work of
$3\cdot K$ with $K=10$ we obtain better accuracy than with $K=640$ over
the whole interval.

In this example
the maximum quadrature errors using the hyperbolas are
smaller than those for the Talbot contours.
Moreover, the hyperbolas allow to choose larger intervals.
On the other hand, the Talbot contours turned out to be
less sensitive to the choice of parameters and the Laplace transform functions
than the hyperbolas.

\subsection{Theoretical error bounds of the contour integral approximations}

For the case of the hyperbola, we obtain in the same way as in Theorem~3 %6.4
of \cite{LoLPS}  the
following error bound  which shows exponential convergence.

\begin{theorem}\label{thm:err}
There are positive constants $C,\,d$, $c_0,\dots,c_4$, and $c$
 such that  at  $t=nh\le T$
the quadrature error in $(\ref{eq:num-int})$ for a hyperbola
$(\ref{Eq.contourparametrisierunghyperbola})$
is bounded by
\begin{eqnarray*}
\|E(\tau,K,h,n)\| &\leq& C \, \dt\,t^{\nu-1}\,(\mu t)^{1-\nu}
\Bigg(  \frac{e^{c_0\mu t}}{e^{2\pi d/\tau}-1}
+ e^{(c_1-c_2\cosh(K\tau))\mu t}
\\
&&
\qquad\qquad \qquad\qquad\
+\ e^{c_3\mu t} \Bigl(1+\frac{c_4\cosh(K\tau)\mu t}{n/2} \Bigr)^{-n/2}\Bigg),
\end{eqnarray*}
if $n \ge c\mu t$ and $\mu t\ge 1$. Here
$\nu$ is the exponent of $(\ref{sector})$.
\end{theorem}

Given an error tolerance $\varepsilon$,
the first term in the error bound becomes
$O(\varepsilon \, \dt \, t^{\nu-1})$ if $\tau$ is chosen so small that
$c_0\mu t - 2\pi d/\tau \le \log \varepsilon$, which requires
an asymptotic proportionality
$
\frac 1\tau \sim \log\frac1\varepsilon + \mu t.
$
For $\mu$ chosen such that
$
\frac{a_1}B \log\frac1\varepsilon \le \mu t \le a_1 \log\frac1\varepsilon
$
with an arbitrary positive constant $a_1$ and with $B> 1$, we obtain that the
second term  is  $O(\varepsilon\, h\,t^{\nu-1})$ if
$c_1 -c_2  \cosh(K\tau) \le -B/a_1$, i.e., with
$
\cosh(K\tau) = a_2
$
for a sufficiently large constant $a_2$. With the above choice of $\tau$,
this yields
$
K\sim \log\frac1\varepsilon.
$
The third term then becomes smaller than $\varepsilon\, \dt\,t^{\nu-1}$ for
$
n \ge c \, \log \frac1\varepsilon
$
with a sufficiently large constant~$c$.
In summary, this gives the following bound
for the required number of quadrature
points on the hyperbola.

\begin{theorem} \label{thm:K}
In $(\ref{eq:num-int})$,
a quadrature error  bounded in norm by $\varepsilon \,\dt\, t^{\nu-1}$
for $nh\in I_\ell$ is obtained
with $K=O(\log\frac1\varepsilon)$. This holds
for $n \ge c \,\log\frac1\varepsilon$ (with some constant $c>0$),
with $K$ independent of $\ell$ and of $n$ and $h$ with $nh\le T$.
\end{theorem}

The approximation is, however, poor for the first few $n$,
as we have seen in the numerical experiments.

We refer to \cite{LoPSch} for an optimized strategy to choose the
parameters $\mu,\tau, K$, which takes also  perturbations in the
evaluations of the Laplace transform into
account.

We expect that a similar result to Theorem~\ref{thm:K} holds also for the
Talbot contours, if the Laplace transform has an
analytic continuation beyond the negative real axis from above and below,
as is the case for the fractional powers considered above.

\section{The fast and oblivious algorithm} \label{Sec.FA}
We now describe the convolution algorithm, concentrating on
Runge-Kutta based convolution quadrature. The algorithm differs
slightly depending on whether we want to compute a convolution or to
solve an integral or integro-differential equation of convolution type.

\subsection{The algorithm for computing convolutions}
The algorithm presented here uses the organisation scheme of
the fast convolution algorithm described in a
step by step manner in~\cite{LuS02}.
A pseudo-code for the algorithm developed in~\cite{LuS02}
can be found in~\cite{HiS04}.

For fixed integer $n\le N$  and a given base $B$ we split the
discrete convolution~(\ref{cq}) or (\ref{rk-cq})
into $L$ sums, where $L$ is the smallest
integer such that $n<2B^{L}$:
\begin{eqnarray*}
u_{n+1} &=&  \sum_{j=0}^{n}  \omega_{n-j}\, g_j =
u_{n+1}^{(0)} + \dots + u_{n+1}^{(L)}
\\
\hbox{with} &&
u_{n+1}^{(0)} = \omega_{0}\, g_{n} \ \mbox{ and } \ u_{n+1}^{(\ell)} =
\sum_{j = b_{\ell}}^{b_{\ell-1}-1} \omega_{n-j}\, g_j
\end{eqnarray*}
for suitable $n=b_0>b_1>\dots>b_{L-1}>b_L=0$.
In view of the approximation intervals~(\ref{intervals}),
the splitting is done in such a way that for fixed $\ell$ in each sum from
$b_{\ell}$ to $b_{\ell-1}-1$, we have $n-j \in [B^{\ell-1},2B^{\ell}-2]$.
The $b_{\ell}=\mathtt{b(\ell)}, \ \ell = 1,\dots,L-1$  are determined
recursively by the  following pseudo-code.
\begin{verbatim}

   L = 1; q = 0;
   for n = 1 to N do
      if 2*B^L == n+1 then L = L+1; endif
      k = 1;
      while mod(n+1,B^k) == 0 & k < L
         q(k) = q(k)+1; k = k+1;
      endwhile
      for k = 1 to L-1 do b(k) = q(k)*B^k; endfor
  endfor

\end{verbatim}
Note that for growing $n$, $b_{\ell}$ is augmented by
$B^{\ell}$ every $B^{\ell}$ steps.
On inserting the integral representation (\ref{W-int}) of the Runge-Kutta
quadrature weights
and the
relation~(\ref{en-rk}), i.e.,
$e_{n-j}(h\lambda)=r(h\lambda)^{n-j}q(h\lambda)$,
we obtain
\begin{eqnarray}
  \label{eq:unell}
 u_{n+1}^{(\ell)} \ &=& \sum_{j=b_\ell}^{b_{\ell-1}-1} \omega_{n-j}\, g_j
  = \sum_{j=b_{\ell}}^{b_{\ell-1}-1} {h\over 2\pi i}
  \int_{\Gamma_{\ell}}   e_{n-j}(h\lambda)\otimes F(\lambda)\, d\lambda\, g_j
\\
\nonumber
  &=&
  {1\over 2\pi i}  \int_{\Gamma_{\ell}}  r(h\lambda)^{n-(b_{\ell-1}-1)}\,
   F(\lambda) \,
     y^{(\ell)}(h\lambda)\, d\lambda
\end{eqnarray}
with
$$
y^{(\ell)}(h\lambda) =
h \sum_{j=b_{\ell}}^{b_{\ell-1}-1}
e_{(b_{\ell-1}-1)-j}(h\lambda) \, g_j\, .
$$
Comparing this formula with (\ref{yrk}), we see that $y^{(\ell)}(h\lambda)$
is the Runge-Kutta approximation to the solution at
$t=b_{\ell-1}h$ of the linear initial-value problem
\begin{equation}\label{ivp-l}
y' = \lambda y + g(t), \quad y(b_\ell h) =0,
\end{equation}
and hence $y^{(\ell)}(h\lambda)$
is computed as such, by Runge-Kutta time-stepping.
The integrals are discretized with the quadrature formula discussed
in Section~\ref{Sec.Contour}:
\begin{equation}\label{u-ell}
u_n^{(\ell)} \doteq  \sum_{k=-K}^K w_k^{(\ell)}
\, r(h\lambda_k^{(\ell)})^{n-b_{\ell-1}+1} \, F(\lambda_k^{(\ell)}) \,
y^{(\ell)} (h\lambda_k^{(\ell)}).
\end{equation}
In the $n$th time step, we thus compute $u_{n+1}^{(\ell)}$ and for subsequent
use we update
the Runge-Kutta solutions to the $(2K+1)L$ initial value problems (\ref{ivp-l})
for the integration points $\lambda_k^{(\ell)}$ on the contours $\Gamma_\ell$
for $\ell=1,\dots, L$, doing one time step from $t_n$ to $t_{n+1}$ in each
of these differential equations.

This algorithm does not keep the history $g_j$ $(j=0,\dots,n)$  in memory. 
For each $\ell=1,\dots,L$ and $k=-K,\dots,K$, it stores
the  Runge-Kutta approximation
to (\ref{ivp-l}) at the current time step, the values
$ w_k^{(\ell)}, \, $ 
$\lambda_k^{(\ell)},\, $
$ r(h\lambda_{k}^{(\ell)}) $,
$ F(\lambda_k^{(\ell)}), \, $ 
$ y^{(\ell)} (h\lambda_k^{(\ell)}), \, $ and
two auxiliary values of the dimension of $y$ needed
for book-keeping purposes (cf.~\cite{LuS02, HiS04}).  
There are
only $(2K+1)L$ evaluations of the Laplace transform $F(s)$.
In the case of real functions
$f(t)$ and $g(t)$ only the real parts of the above sums are needed,
and hence the factor $2K+1$ can be replaced by $K+1$, since the
quadrature points lie symmetric with respect to the real axis.
We recall $L\le \log_B N$ and $K=O(\log\frac1\varepsilon)$, where
$\varepsilon$ is the accuracy requirement in the discretization
of the contour integrals.

In view of the poor approximation of the first convolution quadrature weights
by the discretization of the contour integral, we  evaluate
$u_{n+1}^{(\ell)}$ directly for a few of the first $\ell$, e.g.,
for $\ell=0,1$ with $B=10$. For this we need to keep the $n-b_1+1\le 2B$
values $g_{b_1},\dots,g_n$ in memory, but none of the earlier history
$g_j$ for $j\le n-2B$. We also need the few convolution quadrature weights
$\omega_0,\dots,\omega_{2B-1}$, which may be computed from
(\ref{w-circle-int}) with $2B$ evaluations of the Laplace transform $F(s)$.

For the convolution quadrature based on the
second-order BDF method a similar fast algorithm is obtained by inserting
the formula (\ref{en-bdf2}) for $e_{n-j}(h\lambda)$ in (\ref{eq:unell}).

\subsection{The algorithm for solving integral equations}
\label{subsec:int-eq}
The adaptation of the above algorithm to integral equations such as
\begin{equation}\label{int-eq}
u(t) = a(t) + \int_0^t f(t-\tau)\, g(\tau, u(\tau))\, d\tau\,,\qquad t\ge 0,
\end{equation}
is straightforward for the case of the convolution quadrature based on the
implicit Euler method and the second-order BDF method, which use
solution approximations only on the grid $t=n\dt$. The extension of
the Runge-Kutta based algorithm is, however, less immediate, because
the integral approximation uses the internal stages of the Runge-Kutta
method. Consider a Runge-Kutta based convolution quadrature under the
assumptions of Section~\ref{subsec:rk}.
With the column vector of internal stages
$v_n = (v_{ni})_{i=1}^m$, the discretization of (\ref{int-eq}) reads
\begin{equation}\label{int-eq-rk}
v_n = a_n + \sum_{j=0}^n W_{n-j}\, g_j\, , \qquad n\ge 0,
\end{equation}
with $a_n = \bigl( a(t_n+c_i \dt) \bigr) _{i=1}^m$,
with
weight matrices $W_n$ defined by (\ref{Wn}), and with
$g_j = \bigl( g(t_j+c_i \dt, v_{ji}) \bigr) _{i=1}^m$
depending on the stages $v_{ji}$.
The scheme is implicit in~$v_n$.
The solution at $t_{n+1}$ is  approximated by
the last component of the stage vector~$v_n$,
$$
u_{n+1} = v_{nm}\,.
$$
With the proof of \cite[Theorem 4.1]{LuO93} we obtain that the error of this
approximation over bounded time intervals
is bounded by $O(h^\kappa)$ with $\kappa=\min(p,q+1)$,
where $p$ and $q$ are the classical order
and stage order, respectively, of the underlying
Runge-Kutta method. This estimate holds under the assumption that the
solution is sufficiently smooth.
 It gives orders 3 and 4 for the 2- and 3-stage
Radau IIA methods, respectively. The precise approximation order for the
3-stage method (of classical order 5) may
become larger under appropriate conditions on the nonlinearity and
the convolution kernel, cf.~\cite[Theorem~4.2]{LuO93}.

The weight matrix $W_n$  has the integral
representation, cf.~(\ref{W-int}),
$$
W_n = \frac{\dt}{2\pi i} \int_\Gamma  E_n(h\lambda)\otimes F(\lambda)\,
d\lambda\, ,
$$
where the $m\times m$ matrix $E_n(z)$ is defined by (\ref{En-def}).
By Lemma~2.4 of
\cite{LuO93},  for $n\ge 1$, $E_n(z)$ is the rank-1 matrix given by
\begin{eqnarray*}
E_n(z) &=& r(z)^{n-1} (I-z\A)^{-1} \1 b^T (I-z\A)^{-1}
\\
&=& r(z)^{-1}\, (I-z\A)^{-1}\1 \, e_n(z).
\end{eqnarray*}
These relations permit us to proceed for the history term of
(\ref{int-eq-rk}) as we did for
(\ref{rk-cq}). We split the stage vector $v_n$ as
$$
v_n = a_n + v_n^{(0)} + \dots + v_n^{(L)}
\qquad\hbox{with }\quad v_n^{(\ell)} = \sum_{j=b_\ell}^{b_{\ell-1}-1}
W_{n-j} \, g_j
$$
and obtain, like in (\ref{eq:unell}),
$$
v_n^{(\ell)} = \frac1{2\pi i} \int_{\Gamma_\ell}
r(h\lambda)^{n-b_{\ell-1}}\, (I-h\lambda \A)^{-1}\1\otimes
F(\lambda) \, y^{(\ell)}(h\lambda)  \, d\lambda\, ,
$$
where $y^{(\ell)} (h\lambda)$ is again the  Runge-Kutta approximation
at $t=b_{\ell-1}h$ to
the initial-value
problem (\ref{ivp-l}), now for the inhomogeneity values
$g_j = \bigl( g(t_j+c_i \dt, v_{ji}) \bigr) _{i=1}^m$ in place of
$g_j = \bigl( g(t_j+c_i \dt) \bigr) _{i=1}^m$.
For $\ell\ge 2$ or 3, we thus approximate $v_n^{(\ell)}$ as
$$
v_n^{(\ell)} \doteq  \sum_{k=-K}^K w_k^{(\ell)}
\, r(h\lambda_k^{(\ell)})^{n-b_{\ell-1}} \,
(I-h\lambda_k^{(\ell)}\A)^{-1}\1 \otimes  F(\lambda_k^{(\ell)}) \,
y^{(\ell)} (h\lambda_k^{(\ell)})  \, .
$$
The algorithm stores the same values as before.
The memory requirements for the algorithm are thus independent
of the number of stages $m$ and remain essentially the same
as in the pure convolution case.

\pagebreak[3]

\section{Numerical experiments}
\label{NumEx}
We give two examples to illustrate the  application and behavior
of the fast convolution algorithm.

\subsection{A nonlinear Volterra equation}
We consider a nonlinear Volterra integral equation with weakly singular
kernel from~\cite{Le60},
\begin{equation}
  \label{eq:levinsonNLV}
  u(t) = - \int_{0}^{t}
  \frac{\bigl(u(\tau)-\sin(\tau)\bigr)^3}{\sqrt{\pi(t-\tau)}}\, d\tau\,.
\end{equation}
The convolution quadrature  based on the backward Euler
method gives the implicit discretization
$$
u_{n} = \sum_{j=0}^n \omega_{n-j} \bigl( u_j - \sin(jh) \bigr)^3,
$$
where $\omega_n$ is given by (\ref{omega}) with $F(s)=s^{-1/2}$ and
$\delta(\zeta)=1-\zeta$. To solve the nonlinear equation in each time step
we use  Newton iterations. The history term is computed by the fast
algorithm of the previous section.

We consider also the discretizations based on the
backward differentiation method of order 2, cf.~Section~\ref{subsec:ms}, and
on the 2- and 3-stage RadauIIA implicit Runge-Kutta
methods of orders 3 and 5, respectively; see Sections~\ref{subsec:rk}
and~\ref{subsec:int-eq}.

\begin{figure}[htbp]
  \centering
  \includegraphics[height=0.49\textwidth, angle=270]{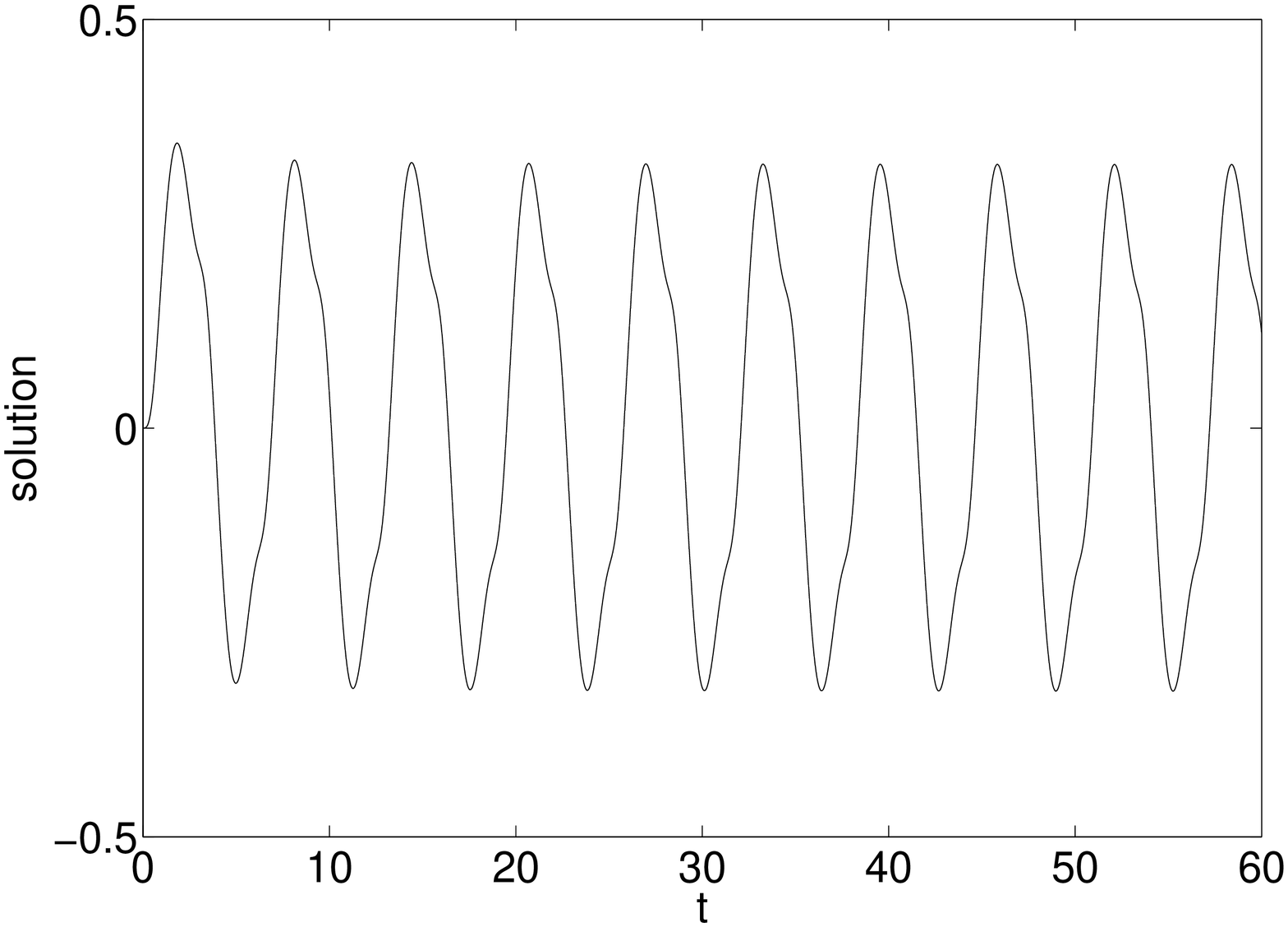}
  \includegraphics[height=0.49\textwidth, angle=270]{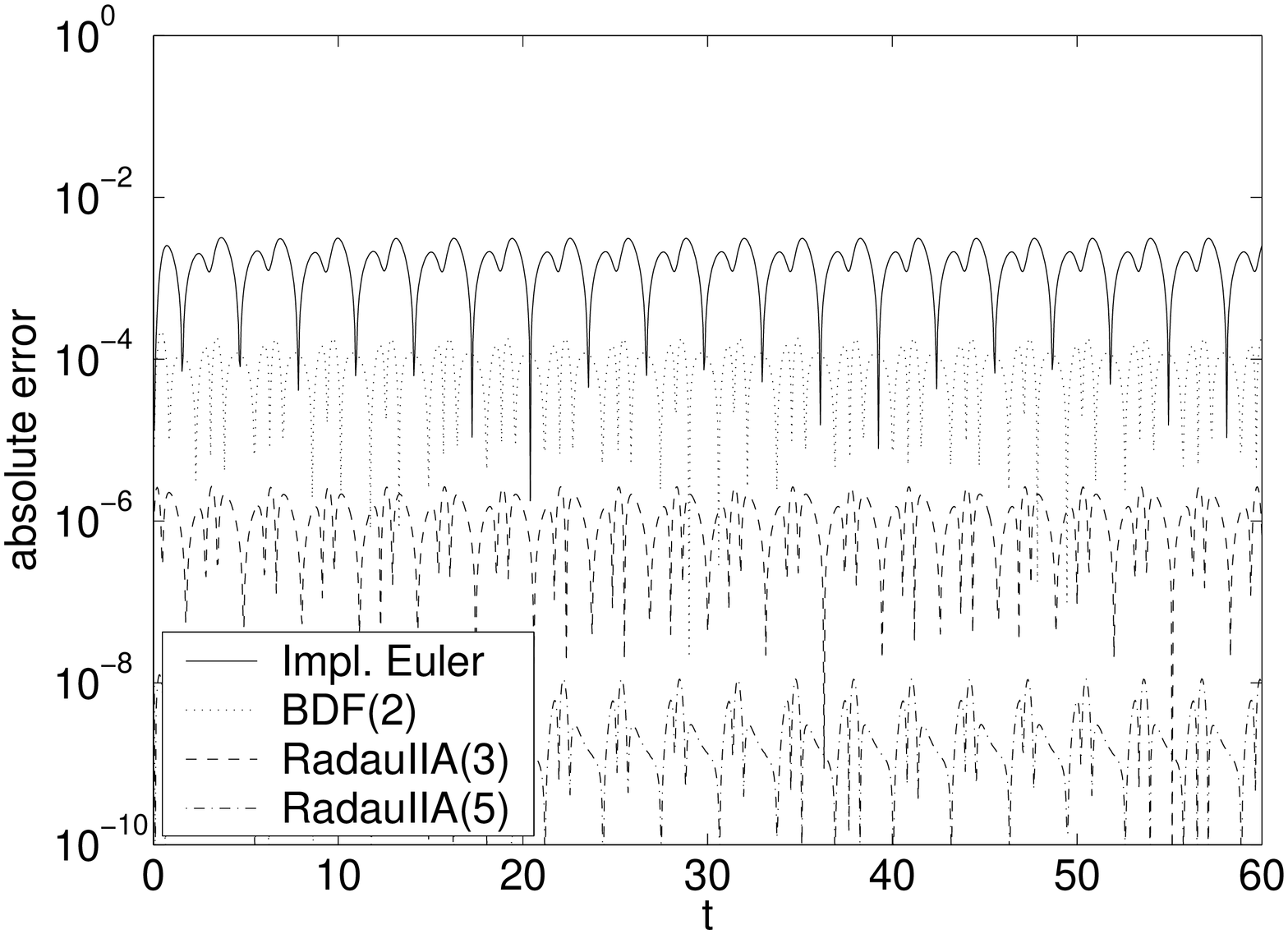}
   \caption{Evolution of the solution over the interval $[0, 60]$ (left) and
   absolute error for different time integration methods
     (right) for $\dt = 0.05$ and $K=30$.}
    \label{fig:errorev}
\end{figure}

In the numerical experiment we use the base $B=5$ and
the Talbot contours with $K = 15$
and $K=30$ and the further parameters as in Section~\ref{subsec:numexp}.
We choose a tolerance of $10^{-12}$ in the Newton method.
The error is calculated with respect to a reference solution,
obtained with $\dt = 0.001$. Figure~\ref{fig:errorev}
shows the evolution of the absolute error and the oscillating solution $u$.

\begin{figure}[htbp]
  \centering
  \includegraphics[height= 0.5\textwidth, width=0.49\textwidth]{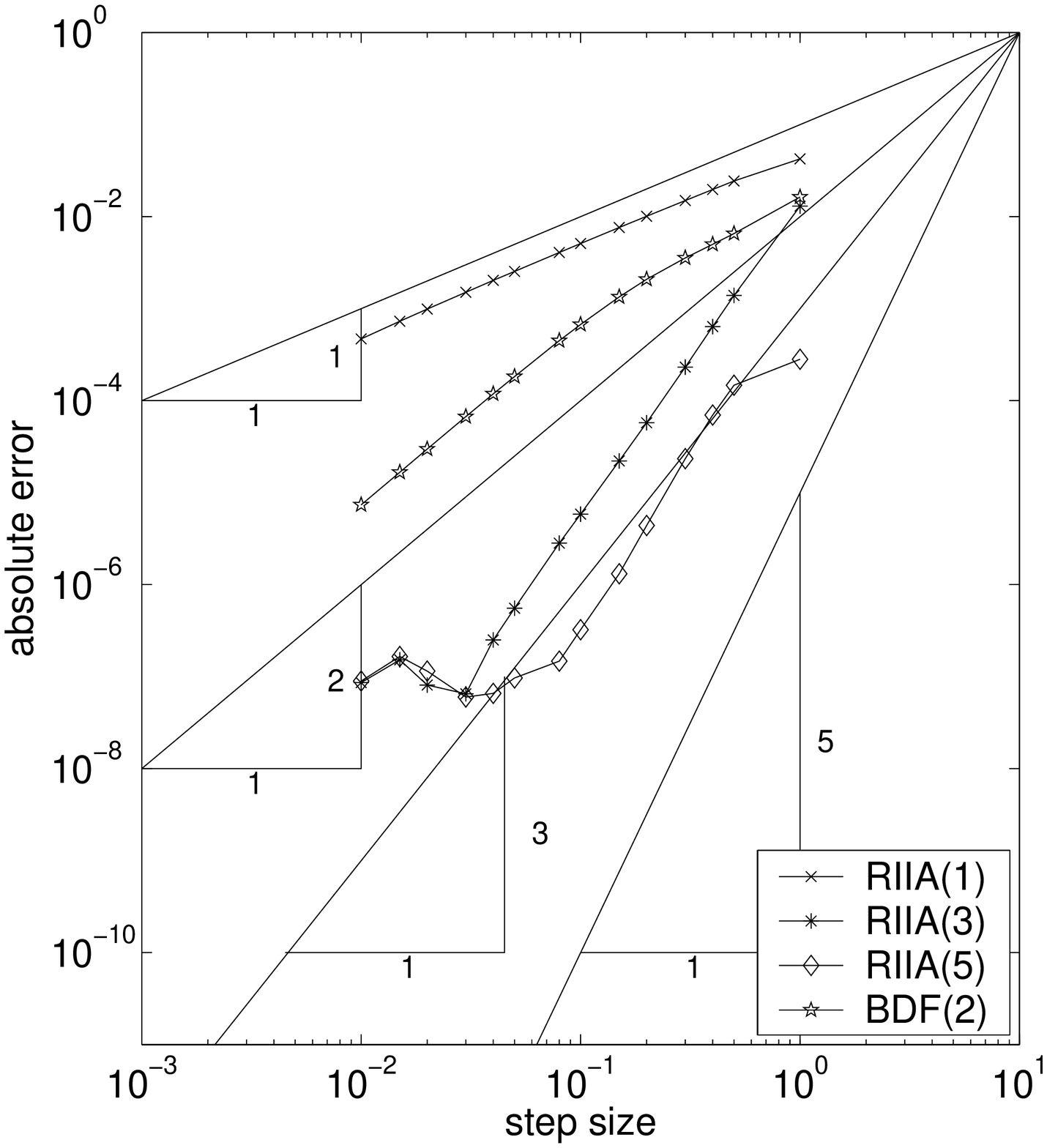}
  \includegraphics[height= 0.5\textwidth, width=0.49\textwidth]{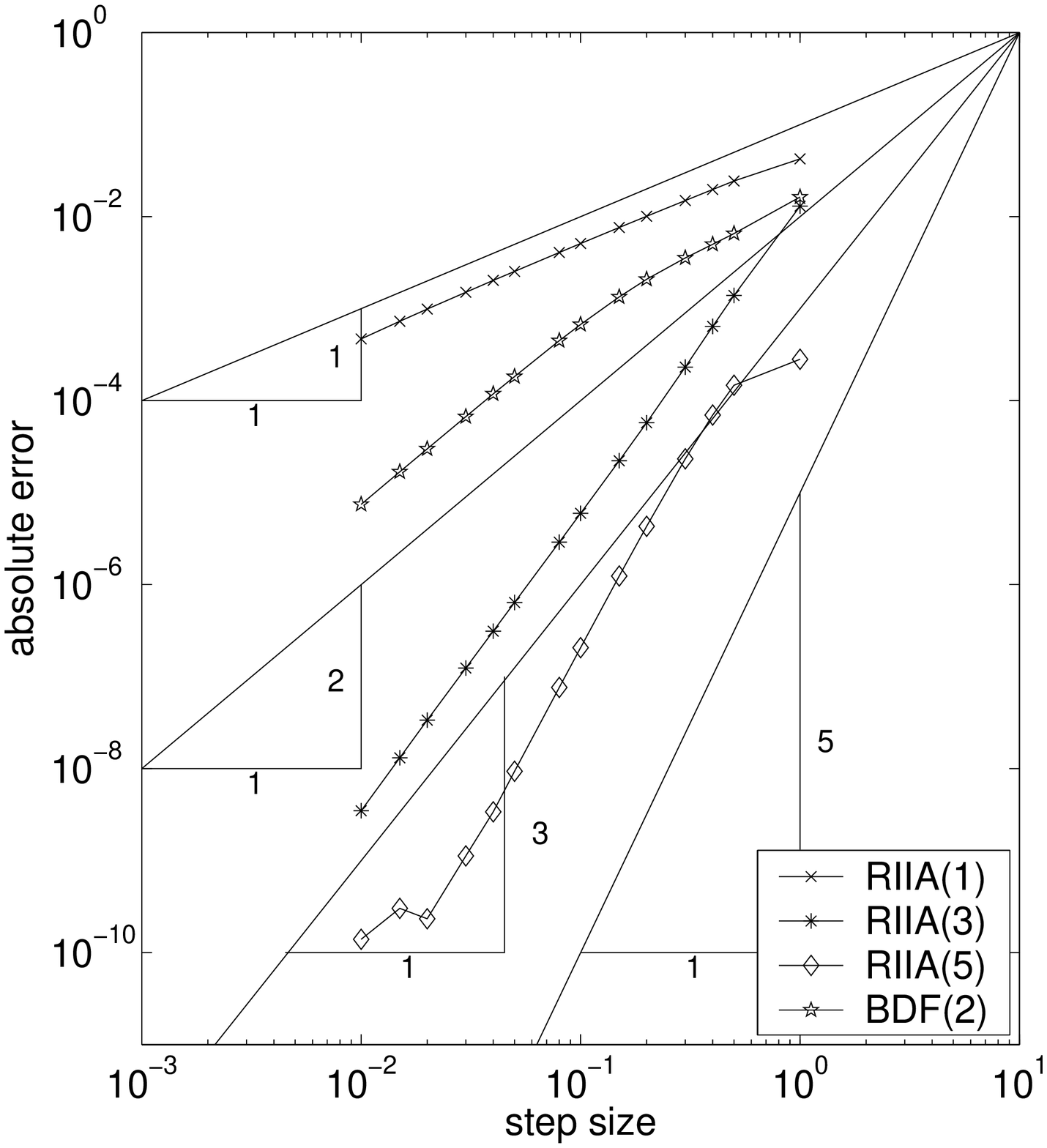}
  \caption{Absolute error vs.~step size $\dt$, for
    different integration methods, with $K=15$ (left) and $K=30$ (right).}
  \label{fig:errordt}
\end{figure}

Figure~\ref{fig:errordt} shows the errors $u_n-u(t_n)$ versus
the step size $\dt$ at time $t_{n}=60$, for $K=15$ and $K=30$.

\begin{figure}[htbp]
  \centering
  \includegraphics[height= 0.5\textwidth, width=0.4\textwidth, angle=270]{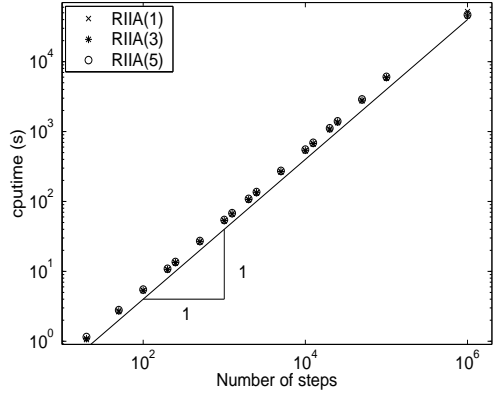}
  \caption{Cpu time in seconds versus the number of integration steps. }
  \label{fig:cputime}
\end{figure}

Figure~\ref{fig:cputime} plots the cpu time
versus the number of integration steps, up to
$10^6$ time steps.
The near-linear growth of the computational work is clearly visible.

\pagebreak[3]

\subsection{Fractional diffusion with transparent boundary conditions}

Here we consider a fractional diffusion equation on the real line;
see, e.g., %Metzler and Klafter
\cite{MeK02} for applications of such
equations in physics and for numerous references.
The equation can be formulated as
\begin{equation}\label{frac-diff}
  u(x,t) - u_0(x) = \int_0^t \frac{(t-\tau)^{\alpha-1}}{\Gamma(\alpha)}
  \, \partial_{xx} u(x,\tau) \, d\tau + g(x,t)
 \quad \mbox{ for } x \in \Rset\, , \ \ t>0
\end{equation}
with the asymptotic condition $u(x,t)\to 0$ for $x\to\pm\infty$,
for an inhomogeneity $g$ with $g(x,0)=0$.
To reduce the computation to a finite domain $x\in [-a,a]$
for initial data $u_0$ and inhomogeneity $g$ with support in $[-a,a]$,
we impose transparent boundary conditions at $x=\pm a$,
which read
\begin{equation}\label{frac-diff-tbc}
u(\pm a, t) = - \int_0^t \frac{(t-\tau)^{\alpha/2-1}}{\Gamma(\alpha/2)}
  \, \partial_{\nu} u(\pm a,\tau) \, d\tau,
\end{equation}
with the outward derivative $\partial_\nu = \pm \partial_x$ at $x=\pm a$.
These boundary conditions are derived with Laplace transform techniques
in the same way as for the wave or the Schr\"odinger equation; see, e.g.,
\cite{H98}.
Space discretisation of~\eqref{frac-diff} is done using second order
finite differences
and a central finite difference to approximate the normal derivative.
With the notation
\[
\delta_{xx} u^{n}_{l} =
\frac{1}{\dx^{2}} \bigl(u^{n}_{l-1} - 2 u^{n}_{l} + u^{n}_{l+1}\bigr)
\ ,\quad \delta_{\nu} u^{n}_{\pm (M-1)} =
\frac{1}{2\dx} \bigl( u^{n}_{\pm M} - u^{n}_{\pm (M-2)}\bigr)
\]
for $a=M\dx$,
the discrete equation approximating~\eqref{frac-diff} is
\begin{gather}
  \label{ueqd}
  \begin{aligned}
  u^{n}_{l} - u^{0}_{l} &=
  \sum_{j=0}^n \omega_{n-j}^{(\alpha)}\,
  \delta_{xx} u^{j}_{l} + g_l^n
  \quad \mbox{ for } l = -(M-1), \dots, M-1 \ ; \ n> 0,
  \\
  u^{n}_{\pm(M-1)} &=
  - \sum_{j=0}^n \omega_{n-j}^{(\alpha/2)}\,
  \delta_{\nu} u^{j}_{\pm(M-1)}~,
  \end{aligned}
\end{gather}
where the  weights $\omega_n^{(\beta)}$
are the convolution quadrature weights for the kernel
$f(t)=t^{\beta-1}/\Gamma(\beta)$ with Laplace transform $F(s)=s^{-\beta}$.

\begin{figure}[htbp]
  \centering
  \includegraphics[height= 0.49\textwidth, width=0.49\textwidth]{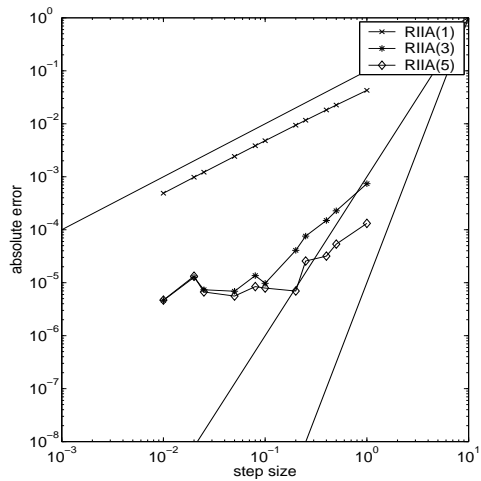}
  \caption{Absolute error vs.~time step, for different integration methods,
  with $K=15$.}
  \label{fig:errorADdt}
\end{figure}

In the numerical example we set $a = 5$ and $M=450$.
We consider the problem with $\alpha = 2/3$ and no inhomogeneity, i.e.,
$g\equiv 0$.
The initial value is $u(x,0) = \exp(-x^{2})$.
Figure~\ref{fig:errorADdt} shows the errors at $t=2$ in dependence on the
step size for the Radau IIA methods of orders 1, 3, 5, obtained with
$B=5$ and $K=15$
in the fast convolution algorithm.
The reference solution is
obtained with the Radau IIA method of order 5, with $\dt = 0.0002$  and $K=40$.
We observe an order reduction for the higher-order
methods, which is due to the temporal non-smoothness of the solution
at $t=0$; cf.~\cite[Sect.\,8]{CuLuPa04}. Nevertheless, the higher-order methods
give much better accuracy.

The work diagram looks almost identical to Figure~\ref{fig:cputime},
showing practically linear dependence of the computational
work on the number of time steps.
The required memory
is less than $200$ entries per spatial grid point for
up to $N\le 10^4$ steps, and less than $300$ entries per grid point
for $N\le 10^6$ steps.
These numbers are halved if we run the algorithm
with $B=10$, $K=10$ instead of $B=5$, $K=15$, as is sufficient for
less stringent accuracy requirements ($\sim 10^{-3}$).

%The solution at $t=2,4,6,8,10$ is shown in Figure~\ref{fig:solAD2o3}
% \begin{figure}[htbp]
%   \centering
%  \includegraphics[width=0.5\textwidth, angle=270]{SolADa0_667NT10000NX801RadauIIA5hyperbola.eps}
%   \caption{Solution  $\alpha = 2/3$.}
%   \label{fig:solAD2o3}
%  \end{figure}

\bigskip\pagebreak[3]

\bibliographystyle{amsplain}

\end{document}